\documentclass[10pt]{article} 
\usepackage{hyperref}
\usepackage{url}

\usepackage[utf8]{inputenc} 
\usepackage{booktabs}
\usepackage{tikz}
\usepackage{amssymb}
\usepackage{amsmath}
\usepackage{graphicx}
\usepackage{stmaryrd}
\usepackage{euscript}
\usepackage{enumitem}
\usepackage[colors,economics]{optsys}
\usepackage{subcaption}

\newcommand{\Bf}{\textsc{Bf}}
\newcommand{\IDS}{\textsc{IDS}}
\newcommand{\FH}{\textsc{FH-Gittins}}
\newcommand{\DeCo}{\textsc{DeCo}}
\newcommand{\Ts}{\textsc{Ts}}
\newcommand{\Klucb}{\textsc{Kl-Ucb}}
\newcommand{\LowerBound}{\textsc{Lb}}
\renewcommand{\multiplier}{\mu}
\newcommand{\job}{a}
\newcommand{\Job}{{\cal A}}
\newcommand{\JOB}{A} 

\newcommand{\good}{\texttt{G}}
\newcommand{\bad}{\texttt{B}}

\newcommand{\Bernoulli}[2]{\mathcal{B}\np{#1,#2}}

\newcommand{\Observation}{Y}

\newcommand{\InformationState}{\pi}

\newcommand{\SIMPLEX}{\Sigma}
\newcommand{\proba}{p}

\newcommand{\prior}{\InformationState}

\newcommand{\REGRET}{{\cal R}}

\usepackage{cancel}

\graphicspath{{Figures/}}
\newcommand{\citep}[1]{\cite{#1}}
\newenvironment{keywords}{\bgroup
  \hsize=\textwidth%
  \noindent\unskip\textbf{Keywords.}\noindent\,\ignorespaces}
{\egroup}

\usepackage{fullpage}

\title{Optimality Gap Analysis of the Decomposition-Coordination\\ Method for Finite Horizon Bandit Problems}

\author{
  Benjamin Heymann\footnote{Criteo Technology, Paris, France}, 
  Michel De Lara\footnote{CERMICS, Ecole des Ponts, Marne-la-Vall\'ee, France},
  Jean-Philippe Chancelier$^*$} 
\date{\today}

\begin{document}

\maketitle

\begin{abstract}
  This study explores multi-armed bandit problems under the premise that the
  decision-maker possesses prior knowledge of the arms' distributions and knows
  the finite time horizon. These conditions render the problems suitable for
  stochastic multistage optimization decomposition techniques.  On the one hand,
  multi-armed bandit algorithms are integral to reinforcement learning and are
  supported by extensive theoretical analysis, particularly regarding regret
  minimization. Meanwhile, stochastic multistage optimization emphasizes
  examining the performance gap between a given policy and the optimal, adapted
  policy. Within this framework, decomposition strategies have been demonstrated
  to efficiently tackle large-scale stochastic multistage optimization
  challenges.  Our empirical findings corroborate the approach's efficacy on
  multi-armed bandit. Most importantly, we contribute a visual interpretation of
  the optimality gap between relaxed and admissible solutions, which lays the
  groundwork for subsequent investigations into the performance of decomposition
  methods.
\end{abstract}


\section{Introduction}

A multi-armed bandit (MAB) is a mathematical model for sequential decision
making under partial feedback. At each round, a decision maker selects an arm,
and then the arm yields a random reward that depends on the intrinsic
characteristics of the arm, which are unknown to the decision maker. The
selection of an arm hence serves two purposes: amassing reward and acquiring
information on the arm, to be used in the future rounds.  For this reason, MABs
embody the well-known exploration-exploitation tradeoff and have concentrated
the attention of the research community for decades.

The first occurrence of MABs in the literature was motivated by clinical
trials~\citep{thompson1933likelihood}, but the rise of the digital economy has
unlocked many new
applications~\citep{NIPS2011_e53a0a29,li2010contextual,weed2016online}.  As
observed in~\citep{chakravorty2014multi}, two schools emerged from the early
works on MABs.  The first school follows~\citep{bellman1956problem} and aims at
maximizing the expected total reward over a discounted horizon, and envisions
the multi-armed bandit as a Markov Decision Process (MDP).  The pioneering
breakthrough is the Gittins Index Theorem~\citep{gittins2011multi} which
provides a way to decompose the problem into tractable subproblems, one per arm.
The second school follows Robbins formulation~\citep{robbins1952some} and seeks
to minimize an expected regret over a finite horizon.  The seminal
work~\citep{lai1985asymptotically,lai1987adaptive} identifies an asymptotically
efficient policy.  Other problem formulations and approaches were proposed
following this
milestone~\citep{lattimore2020bandit,bubeck2012regret,auer2002using}.

In this article, we take the MDP perspective, and aim at maximizing the
intertemporal expected reward or, equivalently, minimizing the Bayesian regret
of a binary bandit over a known, finite timespan.  Theoretically, such a problem
could be addressed with dynamic programming, but this is not feasible because of
the curse of dimensionality: it is well known that the problem size grows
exponentially in the number of arms.  We leverage the ideas
from~\citep{Hawkins:2003,adelman2008relaxations,Carpentier-Chancelier-DeLara-Pacaud:2020,brown2020index}
to show how time decomposition (dynamic programming) can be made compatible with
arm decomposition (Lagrangian relaxation).  This results in a time-dependent
index policy.

Our approach connects to the literature on Lagrangian relaxation of weakly
coupled stochastic dynamic
programs~\citep{Whittle:1988,Hawkins:2003,adelman2008relaxations,brown2020index}.
Our coupling constraint corresponds to the premise that only one arm can be
pulled at any stage.

Asymptotic results for such type of algorithms exist for restless
bandits~\citep{frazier,brown2020index}, when the number of arms that can be
pulled is a constant factor of the time horizon.  However, as far as we know, no
study has been done on MAB with finite horizon. We focus, like
in~\citep{ginebra1999small}, on small time horizons.

We follow~\citep{Carpentier-Chancelier-DeLara-Pacaud:2020}, where it is shown
how a structured, large scale intertemporal maximization problem can be
transformed into a collection of parameterized, simpler, intertemporal
subproblems by relaxing coupling constraints. Thus doing, one obtains a
collection of local value functions, one per subproblems, all functions of a
common coordinating parameter process.  After optimizing this latter, one sums
the local value functions and uses the resulting surrogate global value function
in the online phase of the Bellman equation.  This gives both a theoretical
upper bound for the original maximization problem, and a heuristic online
policy.

This article presents a numerical test bench for the use of the weakly coupled
constraint decomposition method (for short, \DeCo) applied to the Bayesian
bandits problem, and presents an analysis of the optimality gap.  The
decomposition procedure is not new, and it has been in particular used recently
in~\citep{brown2020index} for a setting close to ours (dynamic selection).
Experimentally, what differentiate us from~\citep{brown2020index} is a focus on
a different problem, namely, the Bayesian binary bandit, because (a) this is a
very active area of research for which many solutions have been proposed, (b) it
is a "simpler" problem for which we derive theoretical insights on why those
decomposition methods work so well in practice.  The literature on decomposition
of weakly coupled problems does not tell much on the optimality gap in
nonasymptotic setting. As such, the characterisation of the optimality gap
presented in this paper, while focused on a specific problem, is a step toward
better understanding this class of methods. While it does not fully explain why
\DeCo\ works so well, it is a step into that direction as it provides a visual
interpretation of the gap between the upper bound provided by the relaxed
solution and the performance of \DeCo.

The paper is organized as follows.  In
Sect.~\ref{Decomposition-coordination_method_for_the_bandit_problem}, we
describe the stochastic multi-armed bandit problem, and we show how it can be
treated by decomposition-coordination using a methodology that we refer to as
\DeCo\ .  In
Sect.~\ref{Interpretation_of_the_DeCo_decision_rule_in_term_of_value_of_information_and_reward},
we present an interpretation of \DeCo\ 's decision rule in term of value of
information.  In Sect.~\ref{sec:geometric-representation} we provide a geometric
interpretation of the optimality gap between relaxed solutions and admissible
solutions.We present some algorithmic aspects of \DeCo\ in
Sect.~\ref{sec:deco-algo}.  In Sect.~\ref{sec:numerical-xp}, we show numerically
that \DeCo\ achieves state-of-the-art performance.  It is remarkable that, on
our numerical experiments, \DeCo\ offers performances comparable to
~\citep{russo2014learning}, while also showing empirically good performances for
both small number of arms and for the "many arms" regime.

\section{Preliminaries}

\label{Decomposition-coordination_method_for_the_bandit_problem}

In~\S\ref{Multistage_stochastic_optimal_control_formulation}, we describe the
stochastic multi-armed bandit problem, and show how it can be framed in the
multistage stochastic optimal control formalism.
In~\S\ref{Dynamic_programming_and_arm_decomposition}, we present the arm
decompsition method
in~\citep{brown2020index,Carpentier-Chancelier-DeLara-Pacaud:2020}, and we show
how this control problem can be treated by decomposition-coordination.

\subsection{The Bayesian multi-armed bandit problem}
\label{Multistage_stochastic_optimal_control_formulation}
We now present the (binary) Bayesian multi-armed bandit problem using the
formalism of multistage stochastic optimal control.  For any integers
$r \leq s$, $\ic{r,s}$ denotes the subset $\na{r,r+1,\ldots,s-1,s}$.  We consider a
decision-maker (DM) who selects an arm~$\job$ in a finite set~$\JOB$, at each
discrete time stage~\( t \) in the set \( \ic{0,\horizon{-}1} \), where
$\horizon \geq 1$ is an integer, the \emph{horizon}.  Thus doing, the arm~$\job$
delivers, at the end of the time interval~\( \ClosedIntervalOpen{t}{t{+}1} \), a
random variable\footnote{%
  The shifted index~$t+1$ is here to indicate that the random
  variable~$\va{\Uncertain}_{t+1}^{\job}$ materializes during the time interval
  \( \ClosedIntervalOpen{t}{t+1} \).  }~$\va{\Uncertain}^{\job}_{t+1}$ that
takes two values\footnote{%
  We call these two values ``bad'' (for~\bad), and ``good'' (for~\good), and not
  $\na{0,1}$ to avoid confusion with the possible values for the controls (``do
  not select arm'', ``select arm'').  In fact, we take two values for the sake
  of simplicity, but we could have taken a finite or even infinite number of
  values.  } in the set \( \na{\bad,\good} \) (``bad''~\bad, ``good''~\good) %
and with unknown parameter \(\bar\proba^{\job}_{\good}\), the probability of the
event \( \va{\Uncertain}^{\job}_{t+1}=\good \).  The
parameter~\( \bar\proba^{\job}_{\good} \in [0,1] \) is unknown to the DM, which we
formalize below.

\paragraph{Probabilistic model}

Let \( \SIMPLEX= \bset{ \proba=\np{\proba^{\bad},\proba^{\good}} \in \RR^2}%
  {\proba^{\bad} \geq 0 \eqsepv \proba^{\good} \geq 0 \eqsepv \proba^{\bad} + \proba^{\good} = 1 }
  \) be the one-dimensional simplex\footnote{%
  For the sake of symmetry between outcomes~\bad\ and ~\good,
  we do not identify the simplex~$\SIMPLEX$ with the unit segment~$[0,1]$
  by the mapping \( \SIMPLEX \ni \np{\proba^{\bad},\proba^{\good}} \mapsto
  \proba^{\bad} \in [0,1] \).}. 
For any \( \proba=\np{\proba^{\bad},\proba^{\good}} \in \SIMPLEX \), 
we denote by 
\(
  \Bernoulli{\proba^{\bad}}{\proba^{\good}} = \bigotimes_{t=1}^{\horizon}
  \bp{ \proba^{\bad} \delta_\bad + \proba^{\good} \delta_\good }
\), where $\delta$ denotes a Dirac measure,
the probability on \( \na{\bad,\good}^{\horizon} \)
which corresponds to a sequence of independent (Bernoulli) 
random variables with values in \(  \na{\bad,\good} \). 
For \( \sequence{\proba^{\job}}{\job\in\JOB} =
\sequence{\np{\proba^{\bad\job},\proba^{\good\job}}}{\job\in\JOB}
\in %
\prod_{\job\in\JOB}\SIMPLEX \), we consider the probability
\( \bigotimes_{\job\in\JOB} \Bernoulli{\proba^{\bad\job}}{\proba^{\good\job}} \) on
the product space \( \prod_{\job\in\JOB} \na{\bad,\good}^{\horizon} \), which
corresponds to independence between arms in~$\JOB$.
We denote by \( \EE_{\sequence{\proba^{\job}}{\job\in\JOB}} \) the corresponding mathematical expectation. 
Let $\Delta(\SIMPLEX)$ denote the set of probability distributions
         on the simplex~$\SIMPLEX$. %
We suppose that the DM holds a prior \( \prior^{\job}_{0} \in\Delta(\SIMPLEX) \) over
the unknown~\( \proba^{\job}=\np{\proba^{\bad\job},\proba^{\good\job}} \in \SIMPLEX \), 
for every arm~$\job \in \JOB$.
In practice, we will consider a beta distribution~$\beta\np{n^{\bad},n^{\good}}$
on~$\SIMPLEX$, with positive integers~\( n^{\bad} >0 \) and \( n^{\good}>0\) as
parameters. 

We consider the probability space~$\epro$ where
\(
  \Omega = \prod_{\job\in\JOB}\SIMPLEX \times %
  \na{\bad,\good}^{\horizon} \), 
   \(   \tribu{\NatureField}= \bigotimes_{\job\in\JOB}\borel{\SIMPLEX}
    \otimes    2^{(\na{\bad,\good}^{\horizon})} \) (where \( \borel{\SIMPLEX} \)
    is the Borel $\sigma$-field over the simplex~$\SIMPLEX$),
\(  \PP=\bigotimes_{\job\in\JOB} \prior_0^{\job}\bp{\dd\np{\proba^{\bad\job},\proba^{\good\job}}}
  \otimes %
  \Bernoulli{\proba^{\bad\job}}{\proba^{\good\job}}
\).
Then, \(
\va{\Uncertain}^{\job}=\sequence{\va{\Uncertain}^{\job}_{t}}{t\in\ic{1,\horizon}}
\)  denotes the coordinate mappings for every arm~$\job \in \JOB$,
with \( \va{\Uncertain}^{\job}_{t} \) a random
variable having values in the set \( \na{\bad,\good} \).
For a given family \(
\sequence{\np{\bar\proba^{\job}_\bad,\bar\proba^{\job}_\good}}{\job\in\JOB}
\in \prod_{\job\in\JOB}\SIMPLEX  \)
and for \(
\prior_0^{\job}=\delta_{\np{\bar\proba^{\job}_\bad,\bar\proba^{\job}_\good}} \), 
for every arm~$\job \in \JOB$, the family 
\( \sequence{\va{\Uncertain}^{\job}_{t}}{\job \in \JOB,t\in\ic{1,\horizon}} \)
consists of independent random variables, where \( \va{\Uncertain}^{\job}_{t} \) 
has (Bernouilli) probability
distribution with parameter \( \bar\proba^{\job}_{\good} \in [0,1] \), that is,
\( \PP\bp{\va{\Uncertain}^{\job}_{t}=\bad}=1-\bar\proba^{\job}_{\good} \)
and \( \PP\bp{\va{\Uncertain}^{\job}_{t}=\good}=\bar\proba^{\job}_{\good} \).
With this probabilistic model, we represent the sequential independent outcomes
of \( \cardinal{\JOB} \) independent arms (where \( \cardinal{\JOB} \) denotes the
cardinality of the finite set~$\JOB$).

\paragraph{Decision model}

We consider a sequence \( \va{\Control}=
\sequence{\va{\Control}_{t}}{t\in\ic{0,\horizon-1}} \) 
of random variables (on the probability space~$\epro$), where 
\( \va{\Control}_{t}= \sequence{\va{\Control}_{t}^{\job}}{\job\in\JOB}
  \), \(  \va{\Control}_{t}^{\job} \in \na{0,1}
  \), \(  \forall \job\in\JOB
  \), \(  \forall t\in\ic{0,\horizon-1} \).
Their possible values in~\( \na{0,1} \)
represent that either arm~$\job$ has been selected
at the beginning of the time interval~\( \ClosedIntervalOpen{t}{t{+}1} \) 
($\va{\Control}_{t}^{\job}=1$) or not ($\va{\Control}_{t}^{\job}=0$).
Since, at each given time~$t$, one and only one arm has to be selected, we add the
(almost sure) constraint
\begin{equation}
  \sum_{\job\in\JOB} \va{\Control}_{t}^{\job}=1 
  \eqsepv \forall t\in\ic{0,\horizon-1}
  \eqfinp  
  \label{eq:constraintequality}
\end{equation}
Such way of modelling the selection of a fixed number of arms dates back to the
restless bandit paper~\citep{Whittle:1988}, where Whittle replaces the almost
sure constraint by an expectation.  It has been applied for similar problems
in~\citep{chakravorty2014multi,brown2020index}.

\paragraph{Information and admissible controls}

When the arm~$\job$ has been selected at stage~$t$ (that is, when
$\va{\Control}_{t}^{\job}=1$), the DM observes the outcome, in the set
\( \na{\bad,\good} \), of the random variable~$\va{\Uncertain}_{t+1}^{\job}$.
When the arm~$\job$ has not been selected at stage~$t$ (that is, when
$\va{\Control}_{t}^{\job}=0$), the DM observes nothing (zero~0).
Thus, the DM observes the random variable
\(
\va\Observation_{t+1}=\sequence{\va{\Control}_{t}^{\job}\va{\Uncertain}_{t+1}^{\job}
}{\job\in\JOB}
\), which takes values from the set \( \na{\bad,\good, 0} \)
for all  \( t\in\ic{0,\horizon-1} \).
Then, the admissible controls~\( \va{\Control}=
\sequence{\va{\Control}_{t}}{t\in\ic{0,\horizon-1}} \) %
are those that satisfy
\begin{align}
  \sigma\np{\va{\Control}_{0} }=   \na{\emptyset,\Omega}
  \text{ and }\\
  \sigma\np{\va{\Control}_{t} }
  &\subset
    \sigma\np{\va{\Control}_{0},\va\Observation_{1},\va{\Control}_{1} \ldots,\va{\Control}_{t-1},\va\Observation_{t}}
    \eqsepv \forall t\in\ic{1,\horizon-1}
    \eqfinv   
    \label{eq:filtrations}    
\end{align}
where \( \sigma\np{\va{Z}} \subset \tribu{\NatureField}\) is the $\sigma$-field
generated by the random variable~$\va{Z}$ on the probability
space~$\epro$.   

\paragraph{Random rewards}

We consider given a family
\( \sequence{\InstantaneousCost^{\job}_t}{\job \in \JOB,t\in\ic{0,\horizon-1}} \) of
functions \( \InstantaneousCost^{\job}_t: \na{\bad,\good} \to \RR \), that
represent instantaneous rewards as follows.
When the arm~$\job$ has been selected at stage~$t$ (that is, when
$\va{\Control}_{t}^{\job}=1$), the random
variable~$\va{\Uncertain}_{t+1}^{\job}$ materializes and the DM receives the
payoff
\( 1\times\InstantaneousCost^{\job}_t\np{\va{\Uncertain}_{t+1}^{\job}} =
\va{\Control}_{t}^{\job}\InstantaneousCost^{\job}_t\np{\va{\Uncertain}_{t+1}^{\job}}
\).  When the arm~$\job$ has not been selected at stage~$t$ (that is, when
$\va{\Control}_{t}^{\job}=0$), the DM receives the payoff
\( 0=
\va{\Control}_{t}^{\job}\InstantaneousCost^{\job}_t\np{\va{\Uncertain}_{t+1}^{\job}}
\).
Thus, the total random reward associated with the
\emph{control sequence~\( \va{\Control}= \sequence{\va{\Control}_{t}}{t\in\ic{0,\horizon-1}} \)}
is given by
\(
  \sum_{t=0}^{\horizon-1} \sum_{\job\in\JOB}\va{\Control}_{t}^{\job}
  \InstantaneousCost^{\job}_t\np{\va{\Uncertain}_{t+1}^{\job}}
\). %

\paragraph{State variable}

For any arm $\job\in \JOB$ and
control sequence~\( \va{\Control}= \sequence{\va{\Control}_{t}}{t\in\ic{0,\horizon-1}} \), we set 
\begin{align*}
 \va{N}_t^{\job,\good}=\sum_{s=0}^{t-1}  \va{\Control}_{t}^{\job}
  \1_{\na{\va{\Uncertain}^{\job}_{t+1}=\good}} \mtext{ and }
  \va{N}_t^{\job,\bad}=\sum_{s=0}^{t-1}  \va{\Control}_{t}^{\job}
  \1_{\na{\va{\Uncertain}^{\job}_{t+1}=\bad}}
  \eqfinp
\end{align*}
$ \va{N}_t^{\job,\good}$ (resp.  $\va{N}_t^{\job,\bad}$) represents the
quantities of good (resps. bad) pulls that the decision maker has observed up to
time~$t-1$ (right before making a decision at time~$t$) %
when pulling arm~$\job$, When two different control sequences are envisioned, it
is helpful to be able to explicite the dependency of $\va{N}$ on $\va{U}$, so
that we also introduce the notation
\begin{align}
  \va{N}_t^{\va{U},\job}= \bgp{  \va{N}_t^{\job,\good}
  \eqsepv \va{N}_t^{\job,\bad} }
  \eqfinp
\end{align}

Knowing $\va{N}_t^{\va{U},\job}$, the expected reward of pulling arm~$\job$ can
be deduced by the DM through posterior update. This justifies the introduction
of the following functions \( \ell_t^\job : \NN^2 \to \RR \) given by
\begin{align}
  \forall (n^{\bad\job},n^{\good\job}) \in \NN^2
  \eqsepv
  \ell_t^\job(n^{\bad\job},n^{\good\job})= 
  \frac{{n}^{\bad\job}}{{n}^{\bad\job}+{n}^{\good\job}} \InstantaneousCost^{\job}_t\np{\bad}
  +
  \frac{{n}^{\good\job}}{{n}^{\bad\job}+{n}^{\good\job}}
  \InstantaneousCost^{\job}_t\np{\good}
  \label{def:ellt}
  \eqfinv
\end{align}
which corresponds to the arm expected reward computed according to the updated posterior. 

\paragraph{Optimality criteria in the Bayesian framework}
Let $\Delta(\SIMPLEX)$ denote the set of probability distributions
on the simplex~$\SIMPLEX$. %
We denote by \(
\prior_{0} = \sequence{\prior_{0}^{\job}}{\job \in \JOB}
\in \prod_{\job \in \JOB}\Delta(\SIMPLEX) \)
the family of initial priors, one for each arm, 
and we formulate the following maximization problem ---
where the supremum is taken over \(\va{\Control}=
\sequence{\va{\Control}_{t}^{\job}}{\job\in\JOB, t\in\ic{0,\horizon-1}}
\in \na{0,1}^{\JOB\times\ic{0,\horizon-1}}\), subject to
constraints~\eqref{eq:constraintequality} and \eqref{eq:filtrations},
\begin{subequations}
  \begin{align}
    \VALUE_{0}\np{\prior_{0}}
    =\sup \;
    &
      \int_{\Delta(\SIMPLEX)^{\cardinal{\JOB}}}
      \prod_{\job \in \JOB}\prior^{\job}_{0}\np{\dd\proba^{\job}}
      \EE_{\sequence{\proba^{\job}}{\job\in\JOB}} \bigg[ %
      \sum_{t=0}^{\horizon-1}
      \sum_{\job\in\JOB}\va{\Control}_{t}^{\job}
      \InstantaneousCost^{\job}_t\np{\va{\Uncertain}_{t+1}^{\job}}
      \bigg]
      \label{eq:bandit_problem_criterion}
    \\
    &\text{s.t.}\quad \sum_{\job\in\JOB} \va{\Control}_{t}^{\job}
      =1 
      \eqsepv \forall
      t\in\ic{0,\horizon-1}
      \label{eq:bandit_problem_constraint_as}
    \\
    &\hphantom{\text{s.t.}\quad}
      \sigma\np{\va{\Control}_{t} }
      \subset
      \sigma\np{ \va{\Control}_{0},\va\Observation_{1}, \ldots,\va{\Control}_{t-1},\va{\Observation}_{t}}
      \eqsepv \forall t\in\ic{1,\horizon-1}
      \eqsepv \sigma\np{\va{\Control}_{0} }=   \na{\emptyset,\Omega}
      \eqfinp %
      \label{eq:bandit_problem_constraint_nonanticipatitivity}
  \end{align}
  \label{eq:bandit_problem}
\end{subequations}
We denote by $\mathcal{U}^{ad}$ the set of controls satisfying the
constraints~\eqref{eq:constraintequality} and \eqref{eq:filtrations} (or,
equivalently, constraints~\eqref{eq:bandit_problem_constraint_as} and
\eqref{eq:bandit_problem_constraint_nonanticipatitivity}) of
Problem~\eqref{eq:bandit_problem}.  For any control sequence
$\va{U}\in\mathcal{U}^{ad}$, we set the \emph{intertemporal expected reward} or \emph{total
  reward}
\begin{align}
  \VALUE_{0}^\va{U}\np{\prior_{0}}
  &=
    \int_{\Delta(\SIMPLEX)^{\cardinal{\JOB}}}
    \prod_{\job \in \JOB}\prior^{\job}_{0}\np{\dd\proba^{\job}}
    \EE_{\sequence{\proba^{\job}}{\job\in\JOB}} \bigg[ %
    \sum_{t=0}^{\horizon-1}
    \sum_{\job\in\JOB}\va{\Control}_{t}^{\job}
    \InstantaneousCost^{\job}_t\np{\va{\Uncertain}_{t+1}^{\job}}
    \bigg]
    \eqfinp 
    \label{eq:value_admissible_control}
\end{align}

\subsection{Arm Decomposition}
\label{Dynamic_programming_and_arm_decomposition}

The stochastic optimal control problem~\eqref{eq:bandit_problem} is,
theoretically, solvable by dynamic programming, but
the fact that the computing cost grows exponentially fast with the problem size prevents us from doing so in practice.
One trick leveraged by several
authors~\citep{brown2020index,Carpentier-Chancelier-DeLara-Pacaud:2020} consist
in relaxing the set of constraints~\eqref{eq:bandit_problem_constraint_as}.
This is done in two steps. First the almost sure
constraints~\eqref{eq:bandit_problem_constraint_as} are replaced by constraints
in expectation: $\besp{\sum_{\job\in\JOB} \va{\Control}_{t}^{\job}}=1$.  Second, the
constraints in expectation are replaced by a penalty term in the criterion
$-\mu_t\bp{\besp{\sum_{\job\in\JOB} \va{\Control}_{t}^{\job}}-1}$, where the vector of
penalty parameters $(\mu_t)_{t\in\ic{0,\horizon-1}}$ serve as a lever to push the
system to satisfy the constraint. Those two steps result in the following
modified problem:
\begin{subequations}
  \begin{align}
    \Value^r_0[\mu](\pi_0) =
    \sup \;
    &
      \int_{\Delta(\SIMPLEX)^{\cardinal{\JOB}}} \prod_{\job \in \JOB}\prior^{\job}_{0}\np{\dd\proba^{\job}}
      \EE_{\sequence{\proba^{\job}}{\job\in\JOB}} \bigg[ %
      \sum_{t=0}^{\horizon-1}
      \sum_{\job\in\JOB}\va{\Control}_{t}^{\job}
      \left(  \InstantaneousCost^{\job}_t\np{\va{\Uncertain}_{t+1}^{\job}}-\mu_t \right)
      \bigg] +    \sum_{t=0}^{\horizon{-}1}\multiplier_t 
    \\
    &\text{s.t.}\quad 
      \sigma\np{\va{\Control}_{t} }
      \subset
      \sigma\np{ \va{\Control}_{0},\va\Observation_{1}, \ldots,\va{\Control}_{t-1},\va{\Observation}_{t}}
      \eqsepv \forall t\in\ic{1,\horizon-1}
      \eqsepv \sigma\np{\va{\Control}_{0} }=   \na{\emptyset,\Omega}
      \eqfinp %
  \end{align}
\label{eq:bandit_problem_relaxed}
\end{subequations}
We  refer to $ \Value^r_t[\mu]$ as the \emph{value of the relaxed problem}. 
Because of their Lagrangian interpretation, the vector of penalty parameter $(\mu_t)_{t\in\ic{0,\horizon-1}}$ is often referred to as the vector of multipliers. 
The resulting new problem has the crucial property of being equivalent to  $|A|$ independent, smaller subproblems. As opposed to the initial problem~\eqref{eq:bandit_problem_criterion}, those subproblems can be addressed with dynamic programming. 
More precisely, the value of the relaxed problem writes
\begin{align}
  \label{eq:vr_def}
  \Value^r_0[\multiplier]\bp{\np{n^{\bad\job}}_{\job \in \JOB}, 
  \np{n^{\good\job}}_{\job \in \JOB}}=\\
  \sum_{\job \in \JOB} \Value^{\job}_0\nc{\multiplier}
  \np{n^{\bad\job}_{0}, n^{\good\job}_{0}}
  +
  \sum_{t=0}^{\horizon{-}1}\multiplier_t ,
\end{align}
where $\Value^{\job}_0\nc{\multiplier}
      \np{n^{\bad\job}_{0}, n^{\good\job}_{0}}$ is the value of the subproblem associated with the backward induction:
  \begin{subequations}
  \begin{align}
    \Value^{\job}_\horizon
    \nc{\multiplier}\np{{n}^{\bad\job},{n}^{\good\job}}
    &=
      0      \eqfinv 
      \nonumber
    \\    
    \Value^{\job}_t
    \nc{\multiplier}\np{{n}^{\bad\job},{n}^{\good\job}}
    &=
      \max \Big\{
      \Value^{\job}_{t{+}1}\nc{\multiplier}\np{{n}^{\bad\job},{n}^{\good\job}},
      -\multiplier_t
      + \frac{{n}^{\bad\job}}{{n}^{\bad\job}+{n}^{\good\job}}
      \bp{ \InstantaneousCost^{\job}_t\np{\bad}
      + \Value^{\job}_{t{+}1}\nc{\multiplier}\np{{n}^{\bad\job}+1, {n}^{\good\job}} }
    \\
    \nonumber 
    &\hphantom{=+ max} + \frac{{n}^{\good\job}}{{n}^{\bad\job}+{n}^{\good\job}}
      \bp{ \InstantaneousCost^{\job}_t\np{\good}
      + \Value^{\job}_{t{+}1}\nc{\multiplier}\np{{n}^{\bad\job}, {n}^{\good\job}+1} }
      \Big\}
    \\    
    &=
      \max \Big\{
      \Value^{\job}_{t{+}1}\nc{\multiplier}\np{{n}^{\bad\job},{n}^{\good\job}},
      -\multiplier_t
      + \ell_t^\job(n^{\bad\job},n^{\good\job})
      + \widebarbar{\Value^{\job}_{t{+}1}\nc{\multiplier}}\np{{n}^{\bad\job},{n}^{\good\job}}
      \Big\}
      \label{eq:Pdyna}
      \eqfinv   
  \end{align}
  \label{eq:DP_lambda}
  \end{subequations}
and
where in the last line, we used the notation
\begin{align}
  \forall ({n}^{\bad},{n}^{\good}) \in \NN^2
  \eqsepv 
  \widebarbar{\varphi}({n}^{\bad},{n}^{\good})
  =\frac{{n}^{\bad}}{{n}^{\bad}+{n}^{\good}}
  \varphi \np{{n}^{\bad}+1, {n}^{\good}}
  + \frac{{n}^{\good}}{{n}^{\bad}+{n}^{\good}}
  \varphi\np{{n}^{\bad}, {n}^{\good}+1}
  \eqfinp
  \label{def:barbar}
\end{align}
It is well known~\citep{brown2020index,Carpentier-Chancelier-DeLara-Pacaud:2020}
that $ \Value^r_0[\mu]$ is an upper-bound on the value of the stochastic optimal
control problem~\eqref{eq:bandit_problem} $ \Value_0$, which is what we express
in the next proposition.
\begin{proposition}
  \label{pr:main}
  We have the upper bound
  \begin{align}
    \Value_0
    \bp{\np{n^{\bad\job}_{0}}_{\job \in \JOB}, 
    \np{n^{\good\job}_{0}}_{\job \in \JOB}}
    \leq 
    \inf_{\multiplier  \in \RR^\horizon } \bgp{
    \Value^r_0[\mu]}
    \eqfinp 
    \label{eq:upper_bound_bellman3}  
  \end{align}
  where we identify (by an abuse of notation) 
  \( \Value_0\bp{\np{n^{\bad\job}_{0}}_{\job \in \JOB}, 
    \np{n^{\good\job}_{0}}_{\job \in \JOB}} \)  with 
  the value \( \VALUE_{0}\np{\prior_{0}} \) of problem~\eqref{eq:bandit_problem}
  when the prior \( \prior_{0} = %
  \sequence{\beta\np{n^{\bad\job}_{0},n^{\good\job}_{0}}}{\job\in\JOB}  \) , is a Beta distribution.
\end{proposition}

In what follows we propose to generate a policy using $\Value^r[\multiplier]$
instead of the ``true'' value function in Bellman equation, for a suitable choice of $\mu$.  We call this
procedure \DeCo\ .  This procedure is not new, and has been in particular used
recently in~\citep{brown2020index} for a setting close to ours (dynamic
selection).

It is left as an exercise to check that, when the state of the multi-armed
system is given by
$n= \np{n^{\bad\job}_{t},{n}^{\good\job}_t}_{\job \in \JOB} \in \prod_{\job \in \JOB}
\NN\times\NN$ at time~$t$, the \DeCo\ algorithm selects an arm that maximizes the
value-to-go, or otherwise said, selects an arm $a^\star$ in\footnote{%
  In case of non uniqueness, take any arm in the $\argmax$.
  \label{ft:non-uniqueness}
}
\begin{align}
  \begin{split}
    \forall n \in \NN^{2\cardinal{\JOB}}
    \eqsepv 
    {\JOB^{\sharp}(t,\multiplier,{n})}
    =
    \argmax_{\job \in \JOB}
    \big[ -\Value^{\job}_{t+1} \nc{\multiplier}\np{{n}^{\bad\job},{n}^{\good\job}} 
    +
    \ell_t^\job(n^{\bad\job},n^{\good\job})
    + \widebarbar{\Value^{\job}_{t{+}1}\nc{\multiplier}}\np{{n}^{\bad\job}, {n}^{\good\job}}
    \big]
    \eqfinp       
  \end{split}
  \label{eq:DeCo_policy}    
\end{align}
\section{The value of information interpretation}
\label{Interpretation_of_the_DeCo_decision_rule_in_term_of_value_of_information_and_reward}

Next we provide an interpretation of the online decision rule  of \DeCo\  stated in~\eqref{eq:DeCo_policy}.
If we set 
\begin{align}  
  \delta^{\job}_t\nc{\multiplier}\np{{n}^{\bad\job},{n}^{\good\job}}
  =
  \widebarbar{\Value^{\job}_{t{+}1}\nc{\multiplier}}\np{{n}^{\bad\job}, {n}^{\good\job}}
  -\Value^{\job}_{t+1}\nc{\multiplier}\np{{n}^{\bad\job}, {n}^{\good\job}}
  \label{def:delta}
  \eqfinv
\end{align}
then we can interpret
$\delta^{\job}_t\nc{\multiplier}\np{{n}^{\bad\job},{n}^{\good\job}}$
as the "value of information" in the decomposed subproblem of arm $\job$ when the vector of multipliers  is $\mu$:
it is the incremental performance an optimal policy  gets if given an additional pull outcome. 
The higher $\delta^{\job}_t\nc{\multiplier}\np{{n}^{\bad\job},{n}^{\good\job}}$, the higher a pull of $\job$ increase the expected value one can get from $\job$ in the later rounds in the decomposed problem. 
Hence,  $\delta^{\job}_t\nc{\multiplier}\np{{n}^{\bad\job},{n}^{\good\job}}$  quantify the value of exploration. 
It is insightful to observe that the selected arm  during the \DeCo\ online phase (see Equation~\eqref{eq:DeCo_policy}) is  the one that maximizes
\begin{align}
  \label{eq:ucblike}
  \underbrace{I^\job_t\nc{\multiplier}\np{{n}^{\bad\job},{n}^{\good\job}}}_{\text{index}}=
  \underbrace{\delta^{\job}_t\nc{\multiplier}\np{{n}^{\bad\job},{n}^{\good\job}}}_{\text{value of information (exploration)}}
  +
  \underbrace{\ell^\job_t({n}^{\bad\job},{n}^{\good\job})}_{\text{reward (exploitation)}}
  \eqfinp 
\end{align}
  We recognize an exploration and an exploitation term.
Such exploration  term is reminiscent of the exploration
 term encountered in the Upper Confidence Bound (UCB) algorithms~\cite{auer2002using}. Also,
 \citep{10.1214/aos/1176348788} refers to a learning component in the Gittins index as the
 difference between the index value and the immediate expected reward.
 More recently, the notion of information gain is also important in~\citep{russo2014learning}.
   
 In particular, the definition~\eqref{eq:DP_lambda} of the individual Bellman values in Proposition~\ref{pr:main} becomes
 \begin{align}
   \underbrace{\Value^{\job}_t\nc{\multiplier}\np{{n}^{\bad\job},{n}^{\good\job}}}_{\text{current Bellman value}}
   =& \max \Big\{
      \Value^{\job}_{t{+}1}\nc{\multiplier}\np{{n}^{\bad\job},{n}^{\good\job}},
      -\multiplier_t
      + \ell_t^\job(n^{\bad\job},n^{\good\job})
      + \widebarbar{\Value^{\job}_{t{+}1}\nc{\multiplier}}\np{{n}^{\bad\job},{n}^{\good\job}}
      \Big\}
      \eqfinv  \tag{by~\eqref{eq:Pdyna}}
   \\
   =&
       \underbrace{\Value^{\job}_{t+1}\nc{\multiplier}\np{{n}^{\bad\job},{n}^{\good\job}}}_{\text{future Bellman value at the same state}}
    +\underbrace{\Big(I^\job_t\nc{\multiplier}\np{{n}^{\bad\job},{n}^{\good\job}}    
      - \multiplier_t\Big)^+}_{\text{incremental gain from pulling}}
    \eqfinv    
      \label{eq:value=sum}
 \end{align}
where $x^+=\max(x,0)$,
  which means that the arm is pulled in the decomposed problem only if
  the sum~\eqref{eq:ucblike} of the information gain ($\delta^{\job}_t$) and the
  expected reward ($\ell^\job_t$)
  is greater than $\multiplier_t$. Hence $\multiplier_t$ can be
  interpreted as an 
  equilibrium  price of a ``bandit market''.
  In this bandit market, each bandit is handled by an independent
  profit maximizing agent, and the agent is required to pay the market price~$\multiplier_t$ 
  to pull the arm of her/his bandit at time~$t$.
  This is different but connected to the fair charge metaphore
  proposed in~\citep{10.1214/aoap/1177005588} for the Gittins index.
  The important nuance is that, here, the price depends on a market
  made of several arms whereas, for the Gittins index, the fair charge
  is arm specific. 
\section{Characterization of the optimality gap}
\label{sec:geometric-representation}

In the absence of the constraint~\eqref{eq:constraintequality} of pulling only
one arm, the solutions of the subproblems (one per arm) could be aggregated into an
admissible solution of the original problem~\eqref{eq:bandit_problem}. 
From this perspective, the aggregation of the subproblems solutions constitutes
a solution to the relaxed problem~\eqref{eq:bandit_problem_relaxed}. 

 This section introduces,
 in~\S\ref{Comparative_analysis_of_admissible_and_relaxed_solution_performances},
 a geometric interpretation of the gap between a relaxed solution and any admissible solution (Theorem~\ref{theorem:bound1}).
 We then specialize this result to \DeCo\ (Theorem~\ref{theorem:simu}) in~\S\ref{Optimality_gap_estimate_for_DeCo}.
 The idea is then illustrated with simulation in \S~\ref{sec:illustration}.

\subsection{Comparative analysis of admissible and relaxed solution
  performances}
\label{Comparative_analysis_of_admissible_and_relaxed_solution_performances}

We first show a general result that allows us to compare the total
reward~\eqref{eq:value_admissible_control} generated by an admissible solution
and the total reward generated by a (possibly non admissible) solution to the
relaxed problem~\eqref{eq:bandit_problem_relaxed}.

For $\mu = (\mu_0,\ldots, \mu_{\horizon-1})\in\mathbb{R}_+^{\horizon}$, $\job \in \JOB$, $t\in [0\ldots T-1]$, we define
\begin{align}
  \label{eq:uar}
  \forall \np{{n}^{\bad\job},{n}^{\good\job}}\in \NN^2
  \eqsepv
  u^{a,r}_t\nc{\multiplier}\np{{n}^{\bad\job},{n}^{\good\job}}=
  \begin{cases}
    1 &\text{if }\quad I^\job_t\nc{\multiplier}\np{{n}^{\bad\job},{n}^{\good\job}} - \multiplier_t > 0\eqsepv
    \\
    0 &\text{elsewhere,}%
  \end{cases}
\end{align}

which can be interpreted as a one-arm optimal policy associated with
$\Value^{\job}_t\nc{\multiplier}(n^\job)$, where
$n^\job=\np{{n}^{\bad\job},{n}^{\good\job}}$.  The following result compares the
quantities $\Value^{\va{U}}_0$ (defined in~\eqref{eq:value_admissible_control})
and $\Value^r_0\nc{\multiplier}$(defined
in~\eqref{eq:bandit_problem_relaxed}). Proposition~\ref{pr:main} ensures that
the former is smaller than the latter.  The next result,
Theorem~\ref{theorem:bound1}, quantifies the gap.

\begin{theorem}
  \label{theorem:bound1}
  Let  $\mu = (\mu_0,\ldots, \mu_{\horizon-1})\in\mathbb{R}_+^{\horizon}$,
  and $\va{U}\in\mathcal{U}^{ad}$, and $n\in \mathbb{N}^{2\cardinal{A}}$.
  Then, we have that 
  \begin{align}
    \Value^r_0\nc{\multiplier}(n)-\Value^{\va{U}}_0(n) = 
    \espe_{n}  \Bc{\sum_{t = 0}^{\horizon-1} \sum_{\job \in \JOB }
    \Bp{I^{\job}_t\nc{\multiplier}(\va{N}^{\va{U},a}_t)-\mu_t}\cdot
    \underbrace{\left( u^{a,r}_t\nc{\multiplier}(\va{N}^{\va{U},a}_t)-\va{U}^\job_t\right)}_{\in\{0,-1,+1\}}}
    \label{eq:Bellmandiff}
    \eqfinv
  \end{align}
  where the value function $\Value^r_0\nc{\multiplier}$ is defined in Equation~\eqref{eq:vr_def} and
  the value function $\Value^{\va{U}}_0(n)$ is defined by
  \begin{align}
    \Value^{\va{U}}_0(n) =
    \espe_{n}
    \Bc{\sum_{t = 0}^{\horizon-1} \sum_{\job \in \JOB } \ell^\job_t \np{ \va{N}^{\va{U}}_t} \cdot \va{U}^{\job}_t \np{\va{N}^{\va{U}}_t}}
    \label{def:V0U}
    \eqfinv
  \end{align}
  and $\espe_{n}\nc{\cdot}$ denotes an expectation where the underlying state process takes value $n$ at time $0$,
  that is $\va{N}^{\va{U}}_0=n$.
\end{theorem}

\begin{proof}
  We fix a strategy $\va{U}: \nseqa{n^{\job}}{\job\in \JOB} \to \na{0,1}^{\cardinal{A}}$.
  We start by a two preliminary facts.
  
  \noindent $\bullet$ Using postponed Lemma~\ref{le:def_V0U},
  we obtain that the function $\Value^{\va{U}}_0$ defined in Equation~\eqref{def:V0U}
  coincides with the function $\Value^{\ell,\va{U}}_0$ where the sequence of value functions
  $\nseqp{\Value^{\ell,\va{U}}_{t}}{t\in \ic{0,\horizon}}$ satisfy the Bellman equation~\eqref{eq:Bellman_h},
  that is
  \begin{equation}
    \Value^{\ell,\va{U}}_{t}(n)= {\cal B}_{t+1}^{\ell,\va{U}}[\Value^{\ell,\va{U}}_{t+1}](n)=
    \ell_t\bp{n, {\va{U}}_t(n)} +
    \sum_{\job\in \JOB}
    \widebarbar{\Value^{\va{U}}_{t+1}(n^{(-\job)},\cdot)}(n^\job) \cdot {\va{U}}_t^{\job}(n)
    \eqfinv
    \label{def:BellmanVtU}
  \end{equation}
  where the sequence of mappings $\nseqp{\ell_t}{t\in \ic{0,\horizon-1}}$ is defined, 
  for all $n \in \NN^{2\cardinal{\JOB}}$ and all  $v \in \na{0,1}^{\cardinal{\JOB}}$,  by
  $\ell_t(n,v) = \sum_{\job \in \JOB } \ell_t^{\job}(n^{\job})\cdot v^{\job}$.
    
  $\bullet$ Using again the postponed Lemma~\ref{le:def_V0U}, we obtain that the right
  hand side of Equation~\ref{eq:Bellmandiff} is equal to the value function
  $\Gamma_0^{\va{U}}(n)$ where the sequence of value functions $\nseqp{\Gamma^\va{U}_t}{t\in \ic{0,T}}$ is
  solution of the Bellman equation
  \begin{subequations}
    \label{eq:Bellman_gamma_all}
    \begin{equation}
    \Gamma^\va{U}_T \equiv 0 \text{ and }
    \forall t\in \ic{0,\horizon{-}1}
    \eqsepv
    \Gamma^\va{U}_t = {\cal B}_{t+1}^{\gamma,\va{U}}\nc{\Gamma^{\va{U}}_{t+1}}
    \label{BellmanGamma}
    \eqfinv
  \end{equation}
  where 
  \begin{equation}
    \forall \varphi : \mathbb{N}^{2\cardinal{A}} \to \barRR \eqsepv \forall n \in \mathbb{N}^{2\cardinal{A}}
    \eqsepv 
    {\cal B}_{t+1}^{\gamma,\va{U}}[\varphi](n) =
    \gamma_t\bp{n, \va{U}_t(n)}
    + \sum_{\job\in \JOB}
    \Bp{\widebarbar{\varphi(n^{(-\job)},\cdot)}(n^\job)} \cdot {\va{U}}_t^{\job}(n)
    \eqfinp  \label{eq:Bellman_gamma}
  \end{equation}
  and where the sequence of functions $\nseqp{\gamma_t}{t\in \ic{0,\horizon{-}1}}$ is defined by
  \begin{equation}
    \forall n \in \NN^{2\cardinal{\JOB}}
    \eqsepv
    \forall v \in \na{0,1}^{\cardinal{\JOB}}
    \eqsepv 
    \gamma_t(n,v) = \sum_{\job \in \JOB }
    \bp{I^{\job}_t\nc{\multiplier}(n) -\mu_t}\cdot \np{ u^{a,r}_t\nc{\multiplier}(n) - v^\job}
    \eqfinp
    \label{def:gamma}
  \end{equation}
  \end{subequations}

  Finally, to prove Equation~\ref{eq:Bellmandiff},
  we prove that the sequence of functions $\nseqp{\Value^r_t\nc{\multiplier}-\Value^{\va{U}}_t}{t \in \ic{0,T}}$,
  where the functions $\nseqp{\Value^r_t\nc{\multiplier}}{t \in \ic{0,T-1}}$ are defined by 
  \begin{align}
    \forall t \in \ic{0,T}
    \eqsepv
    \Value^r_t\nc{\multiplier}(n)
    =\sum_{s=t}^{\horizon-1}\mu_s+ \sum_{\job \in\JOB}\Value^{\job}_t\nc{\multiplier}(n^\job)
    \eqfinv
    \label{eq:defVrt}
  \end{align}
  satisfy the Bellman equation~\eqref{eq:Bellman_gamma_all} --- Note that at time
  $0$, Equation~\eqref{eq:defVrt} coincides with Equation~\eqref{eq:vr_def} defining function
  $\Value^r_0\nc{\multiplier}$ ---.

  First, for $t=\horizon$ we immediately check that $\Value^r_{\horizon}\nc{\multiplier}-\Value^{\va{U}}_{\horizon}=0$.
  Second, for $t \in \ic{0,\horizon{-}1}$ and $n \in \NN^{2\cardinal{\JOB}}$, we successively have
  
  \begin{align}
    {\cal B}_{t+1}^{\gamma,\va{U}}
    & \bc{\Value^r_{t+1}\nc{\multiplier}-\Value^{\va{U}}_{t+1}} (n)
      =
      \gamma_t\bp{n, \va{U}_t(n)}
      + \sum_{\job\in \JOB}
      \Bp{\widebarbar{(\Value^r_{t+1}\nc{\multiplier}-\Value^{\va{U}}_{t+1})(n^{(-\job)},\cdot)}(n^\job)} \cdot {\va{U}}_t^{\job}(n)
      \tag{by~\eqref{eq:Bellman_gamma}}
    \\
    &=
      \gamma_t\bp{n, \va{U}_t(n)}
      + \sum_{\job\in \JOB}
      \Bp{\widebarbar{\Value^r_{t+1}\nc{\multiplier}(n^{(-\job)},\cdot)}(n^\job)} \cdot {\va{U}}_t^{\job}(n)
      - \sum_{\job\in \JOB}
      \Bp{\widebarbar{\Value^{\va{U}}_{t+1})(n^{(-\job)},\cdot)}(n^\job)} \cdot {\va{U}}_t^{\job}(n)
    \nonumber \\
    &=
    {\cal B}_{t+1}^{\gamma,\va{U}}\bc{\Value^r_{t+1}\nc{\multiplier}} (n)
      - \sum_{\job\in \JOB}
      \Bp{\widebarbar{\Value^{\va{U}}_{t+1})(n^{(-\job)},\cdot)}(n^\job)} \cdot {\va{U}}_t^{\job}(n)
      \tag{by~\eqref{eq:Bellman_gamma}}
    \\
    &=
      \gamma_t\bp{n, \va{U}_t(n)} +
      \sum_{s=t+1}^{\horizon-1}\mu_s
      +
      \sum_{\job \in \JOB}
      {\cal T}^{\va{U},\job}\bc{\Value_{t+1}^\job\nc{\multiplier}}(n^{\job})
      - \sum_{\job\in \JOB}
      \Bp{\widebarbar{\Value^{\va{U}}_{t+1})(n^{(-\job)},\cdot)}(n^\job)} \cdot {\va{U}}_t^{\job}(n)
      \eqfinv
      \tag{by Lemma~\ref{Bellman_separable} applied to separable function $\Value^r_{t+1}$ given by~\eqref{eq:defVrt} }
      \\
    &=
      \gamma_t\bp{n, \va{U}_t(n)} +
      \sum_{s=t+1}^{\horizon-1}\mu_s
      +
      \sum_{\job \in \JOB}
      {\cal T}^{\va{U},\job}\bc{\Value_{t+1}^\job\nc{\multiplier}}(n^{\job})
      + \sum_{\job\in \JOB} \ell_t^\job\np{n^\job}\cdot {\va{U}}_t^{\job}(n)
      - \Value^{\va{U}}_{t}(n)
      \tag{by~\eqref{def:BellmanVtU}}
    \nonumber \\
    &=
      \gamma_t\bp{n, \va{U}_t(n)}
      + \sum_{s=t+1}^{\horizon-1}\mu_s
      +
      \sum_{\job \in \JOB}
      \Bp{\Value_{t+1}^\job\nc{\multiplier}(n^\job)
      +
      \bp{\widebarbar{\Value_{t+1}^\job\nc{\multiplier}}(n^\job) - {\Value_{t+1}^\job\nc{\multiplier}}(n^\job)
      +\ell_t^\job\np{n^\job } } \cdot {\va{U}}_t^{\job}(n)}
      - \Value^{\va{U}}_{t}(n)
      \tag{by definition of~${\cal T}^{U,\job}$ in~\eqref{def:varphiU1}}
    \\
    &=
      \gamma_t\bp{n, \va{U}_t(n)}
      + \sum_{s=t+1}^{\horizon-1}\mu_s
      +
      \sum_{\job \in \JOB}
      \Bp{\Value_{t+1}^\job\nc{\multiplier}(n^\job)
      +
      I^\job_t\nc{\multiplier}\np{{n}^{\job}} \cdot {\va{U}}_t^{\job}(n)}
      - \Value^{\va{U}}_{t}(n)
      \tag{by~\eqref{def:delta}-\eqref{eq:ucblike}}
    \\
    &=
      \sum_{\job \in \JOB }
      \bp{I^{\job}_t\nc{\multiplier}(n^{\job}) -\mu_t}\cdot \bp{ u^{a,r}_t\nc{\multiplier}(n^{\job}) - {\va{U}}_t^{\job}(n)}
      + \sum_{s=t+1}^{\horizon-1}\mu_s
    \nonumber \\
    &\hspace{3cm}
      +
      \sum_{\job \in \JOB}
      \Bp{\Value_{t+1}^\job\nc{\multiplier}(n^\job)
      +
      I^\job_t\nc{\multiplier}\np{{n}^{\job}} \cdot {\va{U}}_t^{\job}(n)}
      - \Value^{\va{U}}_{t}(n)
      \tag{by definition of $\gamma_t$ in~\eqref{def:gamma}}
    \\
    &=
      \sum_{\job \in \JOB }
      {\Bp{
      \bp{I^{\job}_t\nc{\multiplier}(n^\job) -\mu_t}\cdot { u^{a,r}_t\nc{\multiplier}(n^{\job})}
      + \Value_{t+1}^\job\nc{\multiplier}(n^\job)}}
      + \sum_{s=t+1}^{\horizon-1}\mu_s
      +
      \underbrace{\sum_{\job \in \JOB} \multiplier_t\cdot {\va{U}}_t^{\job}(n)}_{=\mu_t \text{ as } \sum {\va{U}}_t^{\job}=1}
      - \Value^{\va{U}}_{t}(n)
    \nonumber \\
    &=
      \sum_{\job \in \JOB }
      {\Bp{
      \bp{I^{\job}_t\nc{\multiplier}(n^\job) -\mu_t}^{+} 
      + \Value_{t+1}^\job\nc{\multiplier}(n^\job)}}
      + \sum_{s=t}^{\horizon-1}\mu_s
      - \Value^{\va{U}}_{t}(n)
      \tag{by definition of $u^{a,r}_t$ in~\eqref{eq:uar}}
    \\
    &= \sum_{\job \in \JOB }  \Value_{t}^\job\nc{\multiplier}(n^\job)
      +  \sum_{s=t}^{\horizon-1}\mu_s
      - \Value^{\va{U}}_{t}(n)
      \tag{by~\eqref{eq:value=sum}}
    \\
    &= \Value^r_{t}\nc{\multiplier}(n) - \Value^{\va{U}}_{t}(n)
      \eqfinp
      \tag{by definition of $\Value^r_{t}\nc{\multiplier}$ in~\eqref{eq:defVrt}}
  \end{align}
  This ends the proof.
\end{proof}

The following Lemmata~ are instrumental in the proof
of Theorem~\ref{theorem:bound1}.

\begin{lemma}
  \label{le:def_V0U}
  Let be given $\va{U}\in\mathcal{U}^{ad}$ and a sequence of function
  $\nseqp{h_t}{t\in \ic{0,\horizon{-}1}}$ where $h_t:  \NN^{2\cardinal{\JOB}}{\times}\na{0,1}^{\cardinal{\JOB}}\to \barRR$.
  Then, the sequence of value functions $\nseqp{ \Value^{h,\va{U}}_t}{t \in \ic{0,T}}$ defined for $t=\horizon$
  by $ \Value^{h,\va{U}}_T= 0$ and for all $t \in \ic{0,\horizon{-}1}$ by
  \begin{align}
    \Value^{h,\va{U}}_t(n) =
    \espe_{n}
    \Bc{\sum_{s=t}^{\horizon{-}1} h_s \bp{ \va{N}^{\va{U}}_s, \va{U}_s(\va{N}^{\va{U}}_s)}}
    \eqfinv
  \end{align}
  satisfies the Bellman equation
  \begin{subequations}
    \label{eq:Bellman_h}
    \begin{equation}
      \Value^{h,\va{U}}_T \equiv 0 \text{ and }
      \forall t\in \ic{0,\horizon{-}1}
      \eqsepv
      \Value^{h,\va{U}}_t = {\cal B}_{t+1}^{h,\va{U}}\nc{\Value^{h,\va{U}}_{t+1}}
      \eqfinv
    \end{equation}
    where
    \begin{equation}
      \forall \varphi : \mathbb{N}^{2\cardinal{A}} \to \barRR \eqsepv \forall n \in \mathbb{N}^{2\cardinal{A}}
      \eqsepv 
      {\cal B}_{t+1}^{h,\va{U}}[\varphi](n) =  h_t\bp{n, \va{U}_t(n)}
      + \sum_{\job\in \JOB}
      \Bp{\widebarbar{\varphi(n^{(-\job)},\cdot)}(n^\job)} \cdot {\va{U}}_t^{\job}(n)
      \eqfinp  \label{eq:Bellman}
    \end{equation}
  \end{subequations}
\end{lemma}

\begin{proof}
Left to the sagacity of the reader. 
\end{proof}

\begin{lemma}
  \label{Bellman_separable}
  Let $\va{U}\in\mathcal{U}^{ad}$ be given.
  Consider a separable function  $\varphi:  \mathbb{N}^{2\cardinal{A}} \to \barRR$
  defined by $\varphi(n) = \alpha + \sum_{\job \in \JOB} \varphi^{\job}(n^\job)$ then, for all
  $t \in \ic{0,\horizon{-1}}$, we have that
  \begin{equation}
    \forall n \in  \mathbb{N}^{2\cardinal{A}}
    \eqsepv
    {\cal B}^{h,\va{U}}_{t+1}[\varphi](n) =
    h_t\bp{n, \va{U}_t(n)} + \alpha
      +
      \sum_{\job \in \JOB}
      {\cal T}^{\va{U},\job}[\varphi^\job](n^{\job})
      \eqfinv
    \end{equation}
    where the Bellman operator is defined by Equation~\eqref{eq:Bellman} in Lemma~\ref{le:def_V0U} and
    the notation ${\cal T}^{U,\job}$ is defined as follows
    \begin{equation}
      \forall \varphi: \NN^2 \to \barRR \eqsepv
      \forall n \in \NN^{2\cardinal{\JOB}}\eqsepv
      {\cal T}^{\va{U},\job}[\varphi](n) = \varphi(n^{\job}) \bp{1 - \va{U}^{\job}(n)}
      + \widebarbar{\varphi}(n^{\job})\cdot \va{U}^{\job}(n)
      \eqfinp
      \label{def:varphiU1}
    \end{equation}
\end{lemma}
\begin{proof}
  We consider a separable function $\varphi$ given by 
  $\varphi(n) = \alpha + \sum_{\job \in \JOB} \varphi^{\job}(n^\job)$ and we successively compute
  \begin{align}
    {\cal B}^{h,\va{U}}_{t+1}\bc{\varphi} (n)
    &=  h_t\bp{n, \va{U}_{t}(n)}
    + \sum_{\job\in \JOB}
    \Bp{ \widebarbar{\varphi\bp{n^{(-\job)},\cdot}}(n^\job)}
      \cdot {\va{U}}_{t}^{\job}(n^{\job})
      \tag{by~\eqref{eq:Bellman}}\\
    &
      = h_t\bp{n, \va{U}_{t}(n)}
      + \sum_{\job\in \JOB}
      \Bp{\alpha +  \sum_{\job' \in \JOB, \job'\not=\job}
      \varphi^{\job'}(n^{\job'})
      + \widebarbar{\varphi^{\job}}(n^\job)}  \cdot {\va{U}}_{t}^{\job}(n^{\job})
      \tag{as $\varphi$ is separable}
    \\
    & 
      = h_t\bp{n, \va{U}_{t}(n)}
      + \alpha \underbrace{\sum_{\job \in \JOB} {\va{U}}_{t}^{\job}(n^{\job})}_{=1 \text{ as } \va{U}\in\mathcal{U}^{ad}}
      + \sum_{\job\in \JOB}
      \sum_{\job' \in \JOB, \job'\not=\job}
      \varphi^{\job'}(n^{\job'})  \cdot {\va{U}}_{t}^{\job}(n^{\job})
      +  \sum_{\job\in \JOB}
      \widebarbar{\varphi^{\job}}(n^\job)  \cdot {\va{U}}_{t}^{\job}(n^{\job})
    \nonumber \\
    &
      = h_t\bp{n, \va{U}_{t}(n)}
      + \alpha
      + \sum_{\job\in \JOB}
      \varphi^{\job}(n^{\job})  \cdot \bp{1 -  {\va{U}}_{t}^{\job}(n^{\job})}
      +  \sum_{\job\in \JOB}
      \widebarbar{\varphi^{\job}}(n^\job)  \cdot {\va{U}}_{t}^{\job}(n^{\job})
      \tag{by Lemma~\ref{lem:sumAsumA}}
    \\
    & =
      h_t\bp{n, \va{U}_{t}(n)} + \alpha
      +
      \sum_{\job \in \JOB}
      {\cal T}^{\va{U},\job}[\varphi^\job](n^{\job})
      \tag{by~\eqref{def:varphiU1}}
      \eqfinp
  \end{align}
  This ends the proof.
\end{proof}

\begin{lemma}
  \label{lem:sumAsumA}
  Let $g: \JOB \to \RR$ and $j: \JOB \to \RR$ be given and assume that
  $\sum_{\job \in \JOB} j(\job)=1$. Then, we have that 
  \begin{align}
    \sum_{\job \in \JOB} \sum_{\job' \in \JOB, \job' \not=\job} g(\job') j(\job)
    = \sum_{\job \in \JOB}  g(\job)\cdot\bp{1- j(\job)}
    \eqfinp
  \end{align}
\end{lemma}
\begin{proof}
  We successively have
  \begin{align*}
    \sum_{\job \in \JOB} \sum_{\job' \in \JOB, \job' \not=\job} g(\job') j(\job)
    &= \sum_{\job \in \JOB} \sum_{\job' \in \JOB} g(\job') j(\job)\findi{\job\not=\job'} 
      = \sum_{\job' \in \JOB} g(\job') \Bp{\sum_{\job \in \JOB} j(\job)\findi{\job\not=\job'}}
    \\
    &= \sum_{\job' \in \JOB} g(\job') \Bp{\underbrace{\sum_{\job \in \JOB} j(\job)}_{=1} - j(\job')}
    = \sum_{\job' \in \JOB} g(\job') \bp{1 - j (\job')}
    \eqfinp
  \end{align*}
  This ends the proof.
\end{proof}

\subsection{Optimality gap estimate for \DeCo}
\label{Optimality_gap_estimate_for_DeCo}
Next we specialize Theorem~\ref{theorem:bound1} to \DeCo. That is, we consider Theorem~\ref{theorem:bound1}
with the admissible policy used by \DeCo\ denoted by $\text{\DeCo}\nc{\multiplier}$ defined
for all $n=\nseqa{n^{\job}}{\job \in \JOB}\in \NN^{2 \cardinal{A}}$ by (see Equation~\ref{eq:DeCo_policy})
\begin{equation}
  \DeCo\nc{\multiplier}^{\job}(n)
  = \begin{cases}
    1 &\text{if } \job= \job^{\star}(n) \eqfinv \\
    0 &\text{if }  \job\not=  \job^{\star}(n) \eqfinv
  \end{cases}
  \label{eq:decopol}
\end{equation}
where $a^{\star}(n)$ is a unique arm selected in $\argmax_{\job\in\JOB} I_t^\job\nc{\multiplier}(n) $.
\begin{theorem} Let  $\mu = (\mu_0,\ldots, \mu_{\horizon-1})\in\mathbb{R}_+^{\horizon}$,
  and $n\in \mathbb{N}^{2\cardinal{A}}$. Then, we have that 
  \label{theorem:simu}
  \begin{align}
    \Value^r_0\nc{\multiplier}(n) -\Value^{\DeCo\nc{\multiplier}}_0(n)
    &=
      \sum_{t\in [\horizon]}
      \mathbb{E}
      \Bigg[
      \Bp{\mu_t-I^{a^{\star}(\va{N}_{t}^\DeCo)}_t\nc{\multiplier}(\va{N}_{t}^\DeCo)}^+
      +
      \sum_{\job \in A\bp{t,\mu,\va{N}_{t}^{\DeCo,a}}}
      \bp{I_t^\job\nc{\multiplier}(\va{N}_{t}^{\DeCo,a})-\mu_t}
      \Bigg]
      \label{eq:decobell}
      \eqfinv
  \end{align}
  where $A(t,\mu,n)= \nset{a}{a \not={a^{\star}(n)} \wedge I_t^\job\nc{\multiplier}(n)\ge \multiplier_t}$; that is
  the set of arms with index $I_t^{\job}\nc{\multiplier}$ greater than $\mu_t$ which are not selected by
  the \DeCo\ policy. 
\end{theorem}
\begin{proof}
  We first  observe that, by definition of~$u^{a,r}$ (see~\eqref{eq:uar}), we
  have that for all $n=\nseqa{n^{\job}}{\job \in \JOB}\in \NN^{2 \cardinal{A}}$
  \begin{align} 
    \bp{I_t^\job\nc{\multiplier}(n^\job)-\mu_t}
    &\cdot \bp{u^{a,r}_t\nc{\multiplier}(n^\job)-\va{U}_t^\job(n)}
    \nonumber \\
    &=\bp{I_t^\job\nc{\multiplier}(n^\job)-\mu_t}^{+} \cdot
      \bp{1-\va{U}_t^\job(n)}
      + \bp{I_t^\job\nc{\multiplier}(n^\job)-\mu_t}^{-} \cdot
      \va{U}^\job_t(n)
      \eqfinp
      \label{eq:costDeCo}
  \end{align}
  where $x^-=\max(-x,0)$. Now, using the notation $\psi=\DeCo\nc{\multiplier}$, we successively obtain that
  \begin{align}
    \Value^r_0\nc{\multiplier}(n)
    &-\Value^{\psi}_0(n) = 
    \Bgesp{\sum_{t\in [\horizon]} \sum_{\job \in \JOB }\Bp{I^{\job}_t\nc{\multiplier}(\va{N}^{\psi}_t)-\mu_t}
      \cdot \bp{ u^{a,r}_t\nc{\multiplier}(\va{N}^{\psi}_t)- \psi^{\job}_t\nc{\multiplier}
      (\va{N}^{\psi}_t)}}
    \tag{by~\eqref{eq:Bellmandiff}}
    \\
    &= 
      \Bgesp{\sum_{t\in [\horizon]} \sum_{\job \in \JOB }
      \bp{I^{\job}_t\nc{\multiplier}(\va{N}^{\psi}_t)-\mu_t}^{+}
      \cdot \bp{1 - \psi^\job_t (\va{N}^{\psi}_t)}
      + \bp{I^{\job}_t\nc{\multiplier}(\va{N}^{\psi}_t)-\mu_t}^{-}
      \cdot
      \psi^\job_t(\va{N}^{\psi}_t)}
      \tag{by~\eqref{eq:costDeCo}}
    \\
    &= 
      \Bgesp{\sum_{t\in [\horizon]}
      \Bp{
      \bp{I^{\job^{\star}(\va{N}^{\psi}_t)}_t\nc{\multiplier}(\va{N}^{\psi}_t)-\mu_t}^{-}
      +
       \sum_{\job \in \JOB,  \job\not=a^{\star}(\va{N}^{\psi}_t)}
      \bp{I^{\job}_t\nc{\multiplier}(\va{N}^{\psi}_t)-\mu_t}^{+}}}
      \tag{by~\eqref{eq:decopol}}
      \eqfinv
  \end{align}
  which immediately gives Equation~\eqref{eq:decobell} and ends the proof.
\end{proof}
  
Since $ \Value^r_0\nc{\multiplier}(n)$ is an upper bound for the performance of
any admissible control, Theorem~\ref{theorem:simu} indicates that the optimality
gap is upper bounded by the sum of two terms.  One term corresponds to the
situations encountered by \DeCo, when no arms has an index greater than $\mu_t$.
The other term corresponds to the situations when strictly more than one arm has
an index greater than $\mu_t$.

\subsection{Illustration}
\label{sec:illustration}

\begin{figure}[hbtp]
  \begin{center}
    \begin{subfigure}[b]{0.4\textwidth}
      \caption{run 1}
      \includegraphics[width=\textwidth]{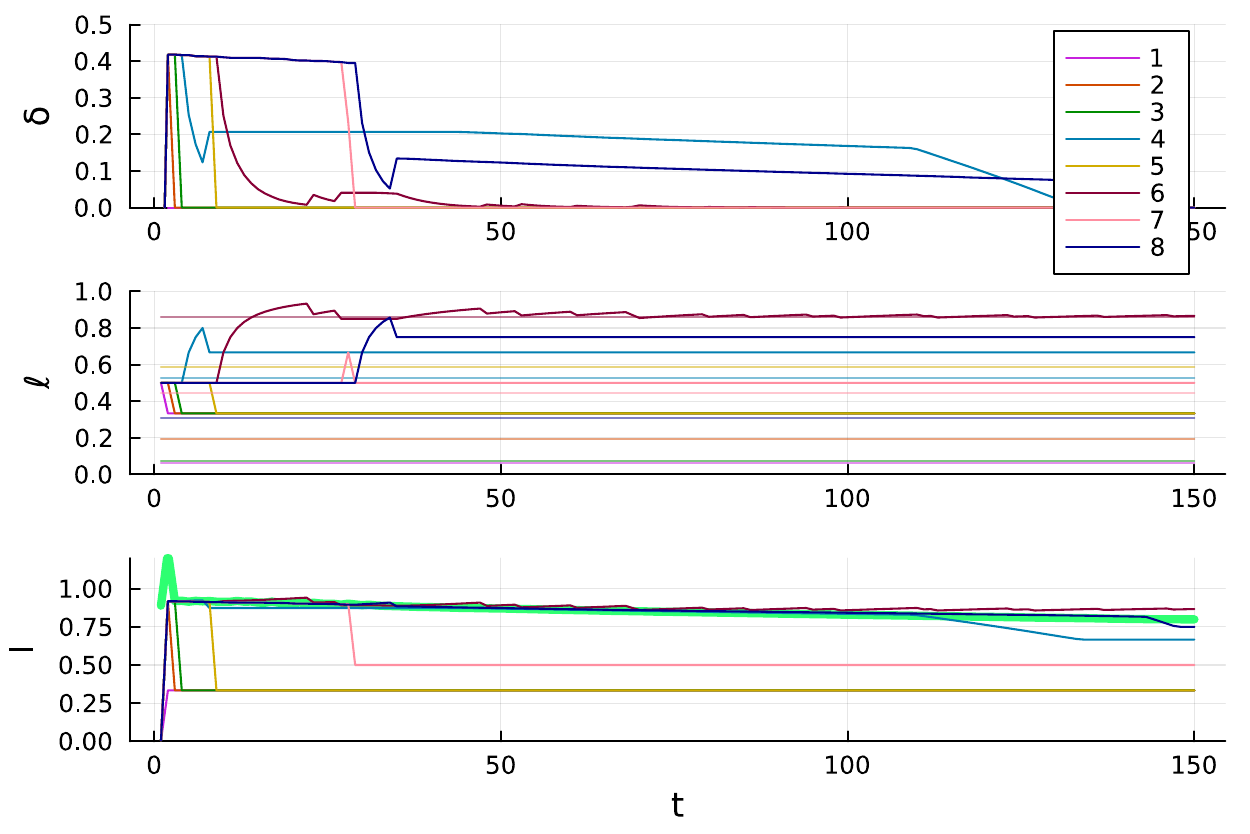}%
    \end{subfigure}%
    \hfill%
    \begin{subfigure}[b]{0.4\textwidth}
      \caption{run 2}
      \includegraphics[width=\textwidth]{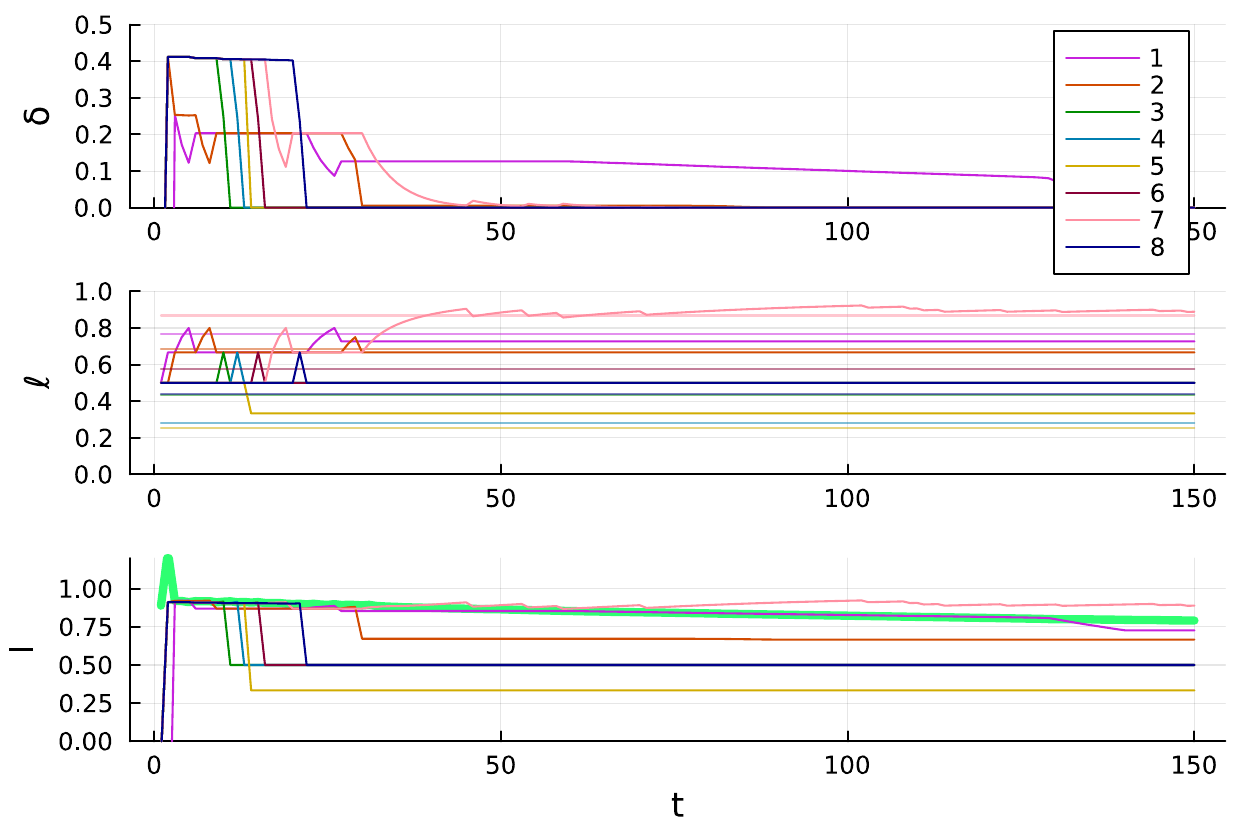}%
    \end{subfigure}
    
    \begin{subfigure}[b]{0.4\textwidth}
      \caption{run 3}
      \includegraphics[width=\textwidth]{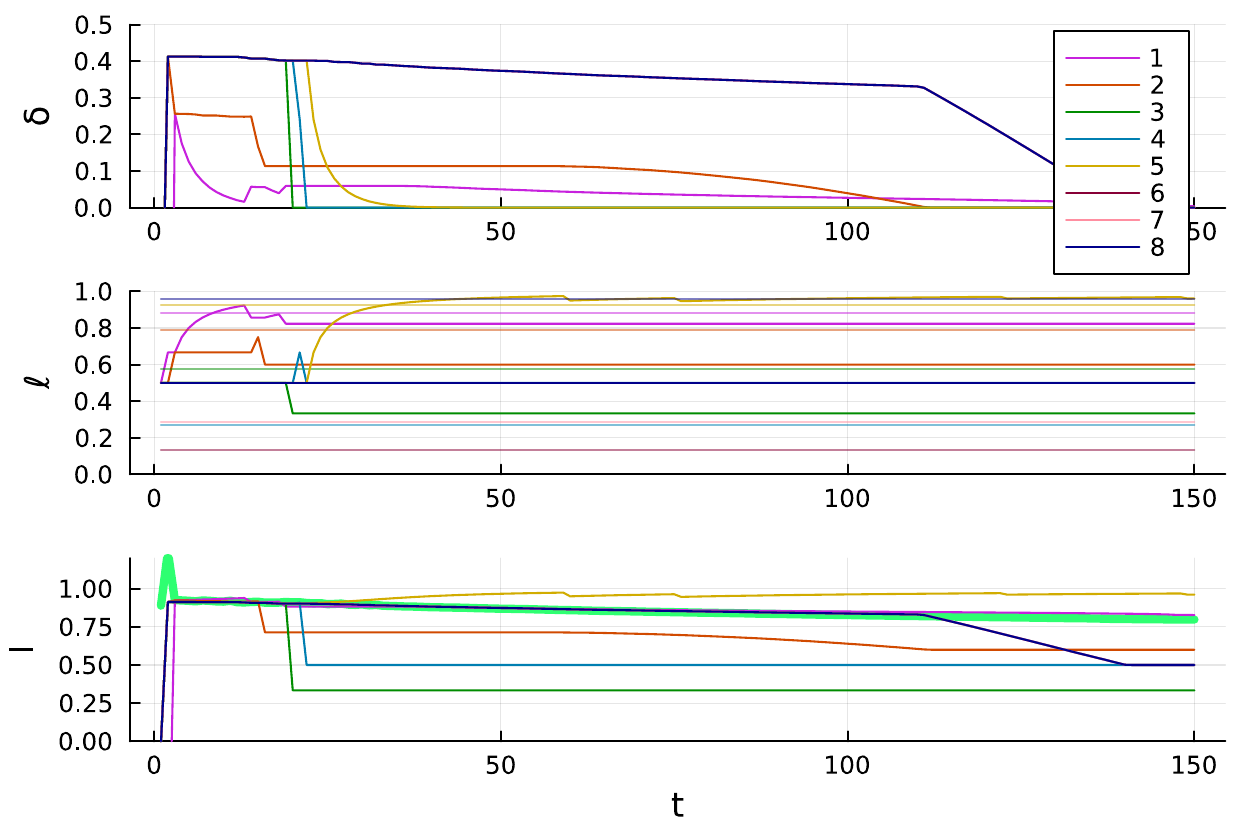}%
    \end{subfigure}%
    \hfill%
    \begin{subfigure}[b]{0.4\textwidth}
      \caption{run 4}
      \includegraphics[width=\textwidth]{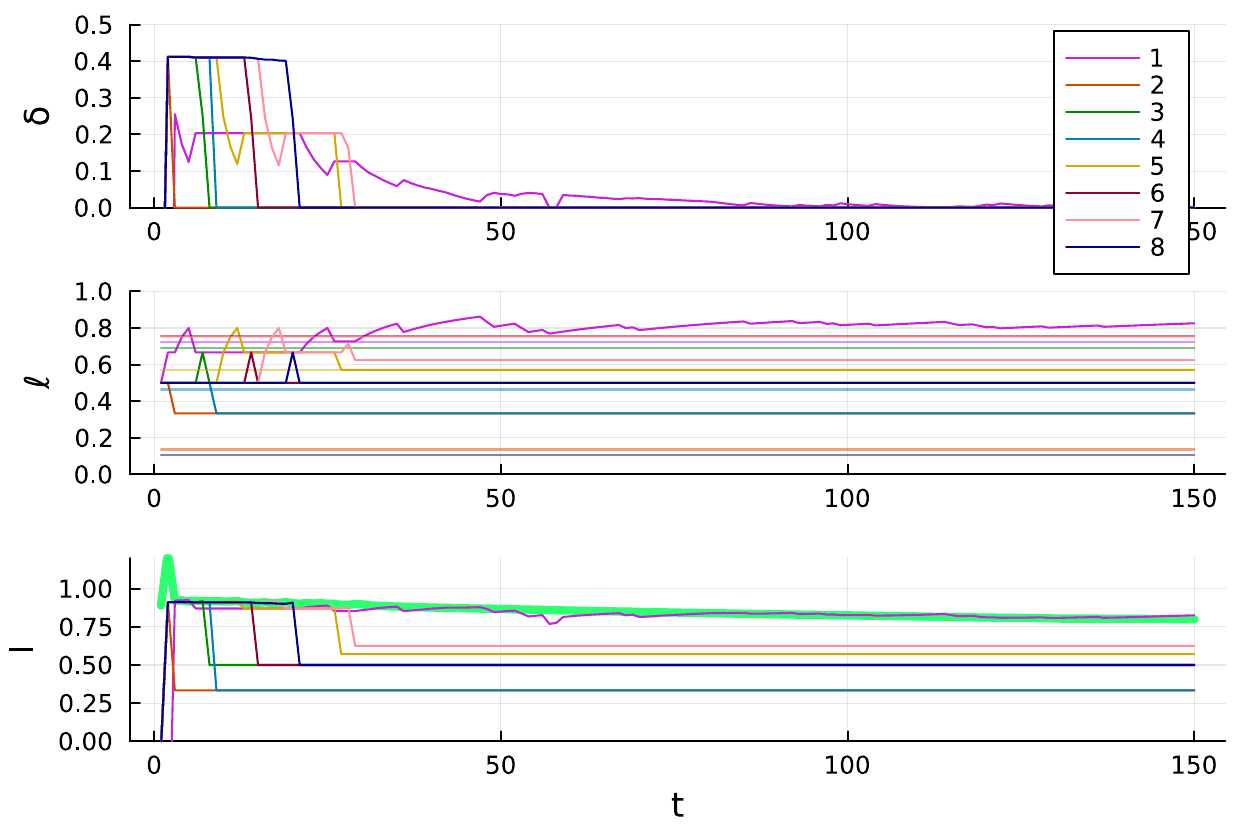}
    \end{subfigure}%
  \end{center}
  \caption{Four simulations of trajectories generated with the controls obtained from \DeCo\ for a uniform prior and 8 arms. 
  In each of the four subfigures, we display the value of information~$\delta$
  (top), the expected reward~$\ell$ (middle) and the index~$I$ (bottom) obtained with \DeCo.
  The green line corresponds to the value of the
  component~$\mu_t$ of the multiplier at the end of the dual gradient procedure.}
\label{fig:result-main2}
\end{figure}

We illustrate Theorem~\ref{theorem:simu} with a  numerical simulation reproduced four times. We use \DeCo\ to produce an approximation of the optimal Bayesian control when 8 arms are uniformly sampled from $[0,1]$, then  we sample uniformly the values of 8 arms in $[0,1]$, and run a simulation. 
We show the results  in Figure~\ref{fig:result-main2}.
The optimality loss can be ``read'' directly from the plot. 

\section{Algorithmic aspects of \DeCo\ }
\label{sec:deco-algo}
The \DeCo\ algorithm is made of an offline computation phase (described in~\S\ref{Offline_phase_of_the_DeCo_algorithm})
and of an online computation phase (described
in~\S\ref{Online_phase_of_the_DeCo_algorithm})
phases as follows.
The offline phase is summarized in Figure~\ref{decompo-algo}: it
consists in the minimization of a dual function~$\varphi$, where each evaluation
of~$\varphi$ relies on solving $\cardinal{\JOB}$~independent Bellman equations. The 
minimization  in~\eqref{eq:upper_bound_bellman3} can be performed by gradient descent.
Then, the online phase consists in using the upper bound
function~\eqref{eq:upper_bound_bellman3} as a proxy for the Bellman value, while
ensuring that only one arm is pulled at each step.
In~\S\ref{Computational_complexity}, we discuss the computation cost of \DeCo\ .

\subsection{Offline phase of the \DeCo\ algorithm}
\label{Offline_phase_of_the_DeCo_algorithm}

The offline phase of the \DeCo\ algorithm is %
the minimization of the upper bound $\Value^r_0[\mu]$ with respect to $\mu$ (see \eqref{eq:upper_bound_bellman3}), 
for a family 
\( \prior_{0}=\sequence{\prior_{0}^{\job}}{\job\in \JOB}=
\sequence{\beta\np{n^{\bad\job}_{0},{n}^{\good\job}_{0}}}{\job \in \JOB} 
 = \beta\np{n_{0}} \)
of beta priors. 
\begin{figure}[hbtp]
  \begin{center}
    {\includegraphics[width=0.5\textwidth]{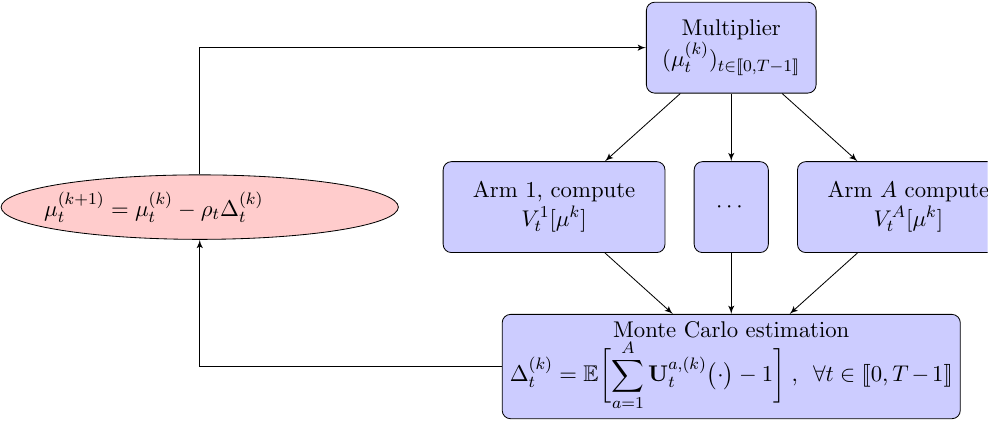}}%
    \caption{The decomposition coordination algorithm (\DeCo)}
    \label{decompo-algo}
  \end{center}
\end{figure}
The algorithm is summarized in Figure~\ref{decompo-algo} and its four 
steps are as follows. 
\begin{enumerate}[label=$(\text{S}_{\arabic*})$, ref=$\text{S}_{\arabic*}$,topsep=\parsep,itemsep=0ex]
\item
  Choose an initial vector~$\multiplier^{(0)} \in \RR^\horizon$ of multipliers.
\item 
\label{step2} 
At iteration $k$, given a vector~$\multiplier^{(k)} \in \RR^\horizon$ of multipliers, compute the collection
  $\sequence{\Value^{\job}_t\nc{\multiplier^{(k)}}}{t\in\ic{0,\horizon}, \job \in \JOB}$.
  The computation is performed in parallel, arm per arm. Note that $\Value^{\job}_t\nc{\multiplier^{(k)}}$, the Bellman value function at time $t \in \ic{0,\horizon}$,
  is to be evaluated only on the finite grid $\nset{\np{n^{\bad\job}_{0}+n^{\bad\job},{n}^{\good\job}_{0}+{n}^{\good\job}}}
  {{n}^{\bad\job}+{n}^{\good\job}\le t}$.
  Note also that, when all the arms share the same prior and the same
  instantaneous reward, a unique sequence of Bellman value functions is to be computed, that
  is, all the arms share the same sequence of Bellman value functions. 
\item
  Once gotten the collection of  $\sequence{\Value^{\job}_0\nc{\multiplier^{(k)}}}{\job \in \JOB}$
  of decomposed Bellman value functions at iteration~$k$,
  update the vector of multipliers by a gradient step to obtain~$\multiplier^{(k+1)}$. The gradient of the dual function
  $\varphi$ with respect to the multipliers is obtained by computing the expectation of the dualized constraint as
  formulated in Problem~\eqref{eq:upper_bound_bellman3}
  (see~\citep{Carpentier-et-all-2018} for more details). Numerically, the expectation is obtained
  by Monte Carlo simulations. \footnote{The gradient phase to minimize~\eqref{eq:upper_bound_bellman3} can be replaced by a more sophisticated algorithm
  such as the conjugate gradient or the quasi-Newton method. }
In some of our numerical experiments, we use  a solver (limited memory \textsc{Bfgs}) of the \textsc{Modulopt}
  library from \textsc{Inria}~\citep{gilbert2007libopt}. To obtain a
  global $O(T^3)$  running
  time, the computing  budget 
  allocated to this iterated  gradient phase does not depend on~$\horizon$.
  
\item Stop the iterations (stopping criterion) or go back to step~\ref{step2} with multiplier $\multiplier^{(k+1)}$.
\end{enumerate}

\subsection{Online phase of the \DeCo\ algorithm}
\label{Online_phase_of_the_DeCo_algorithm}

In the offline phase, we obtain a multiplier $\mu$ and a collection of value function $V^\job$.
During the online phase, the arms are selected according to the rule specified in Equation~\eqref{eq:DeCo_policy}.

The structure of policy~\eqref{eq:DeCo_policy} is that of a
\emph{nonstationary\footnote{%
If we had considered an infinite horizon, we would have obtained a
(stationary) index policy.} index policy}.
Indeed, the right hand side in~\eqref{eq:DeCo_policy} is a quantity that depends
only on~$t$ and on the state \( \np{{n}^{\bad\job}_t,{n}^{\good\job}_t} \) of arm~$\job$ at
time~$t$. 
The \DeCo\ policy used in numerical experiments is the
policy~\( \Job\opt\nc{\multiplier\opt} \) in~\eqref{eq:DeCo_policy},
where $\multiplier\opt$ is given by the offline phase
(described in~\S\ref{Offline_phase_of_the_DeCo_algorithm})
of the \DeCo\ algorithm.

\subsection{Computational complexity}
\label{Computational_complexity}
  
Solving the maximization problem~\eqref{eq:bandit_problem}, that is, computing $\VALUE_{0}\np{\prior_{0}}$
for a given prior (like, for instance, the uniform distribution given by the beta distribution~$\beta(1,1)$ for all arms)
can be done using Dynamic programming.
This is however, only possible for relatively small instances of problem~\eqref{eq:bandit_problem},
that is, for a limited number~$\cardinal{\JOB}$ of arms
and a limited time horizon~$\horizon$. We recall here that solving the problem for
$\cardinal{\JOB}$ arms requires solving a Bellman equation with a state of
dimension $2\cardinal{\JOB}$ (a state described by two integers per arm), which
implies an exponential increase~$O((2|A|)^\horizon)$ in computational cost with respect
to~$\cardinal{\JOB}$. This is an instance of what Richard Bellman
referred to as the \textit{curse of dimensionality}.
For~\FH\ , ~\citep{nino2011computing} provides methodologies to
  compute FH-Gittins indices that are $O(\horizon^6)$ in time complexity.
 
  By contrast, for \DeCo, each dynamic programming phase costs~$O(\horizon^3)$ in running
  time: indeed for each time $t=1,\ldots, \horizon$,  we need a grid of $\horizon\times \horizon$ for the 2 dimensional prior
  parameter that counts the number of  successes and failures, hence
  $O(\horizon^3)$ in total, hence the 
  complexity is $O(K \horizon^3)$, where $K$ is the number of
 arms with different parameters.
  
\section{Numerical experiments}
\label{sec:numerical-xp}

  In this section, we report numerical experiments for some short
  horizons instances of the  Bayesian Bernoulli MAB problem.
For the sake of reproducibility, we have performed two separate 
implementations with two different languages, one in
Julia~\citep{bezanson2017julia}, the other in Nsp~\citep{nsp}.
On the instances we have tested, \DeCo\ competes with state-of-the-art policies for MABs.

First, in Table~\ref{BFvsDeco}, we compare the performance of~\DeCo\ against the
brute force approach~\Bf\ and ~\FH, 
(because of the computation cost of \Bf , \FH\ is used as a proxy
supposed to be close to the optimal solution).
We observe that the performance of~\DeCo\ is close\footnote{It might happen that~\DeCo\  empirical
 average is better than~\Bf, but this is due to the simulation noise.} to the optimal
solution while keeping the computational cost reasonable (at most 1.3~second).

\begin{figure}[hbtp]
  \begin{center}
    \begin{subfigure}[b]{0.4\textwidth}
      \caption{$2$ arms}
      \includegraphics[width=\textwidth]{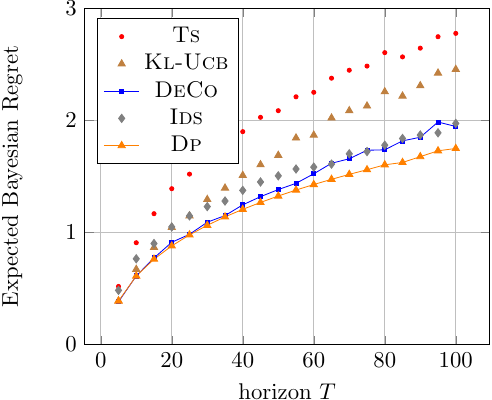}%
    \end{subfigure}%
    \hfill%
    \begin{subfigure}[b]{0.4\textwidth}
      \caption{$5$ arms}
      \includegraphics[width=\textwidth]{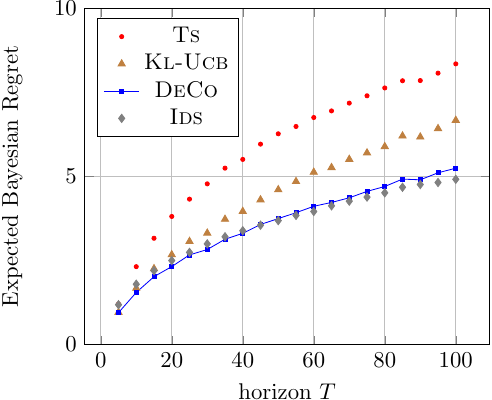}%
    \end{subfigure}
    
    \begin{subfigure}[b]{0.4\textwidth}
      \caption{$10$ arms}
      \includegraphics[width=\textwidth]{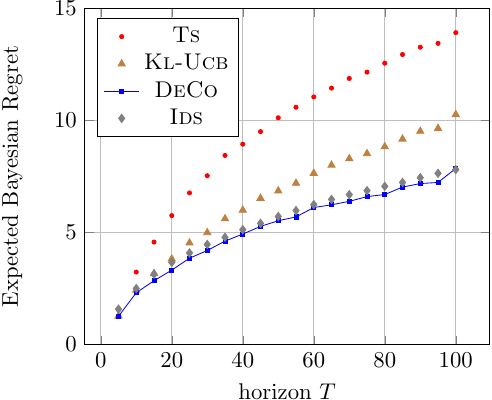}%
    \end{subfigure}%
    \hfill%
    \begin{subfigure}[b]{0.4\textwidth}
      \caption{$20$ arms}
      \includegraphics[width=\textwidth]{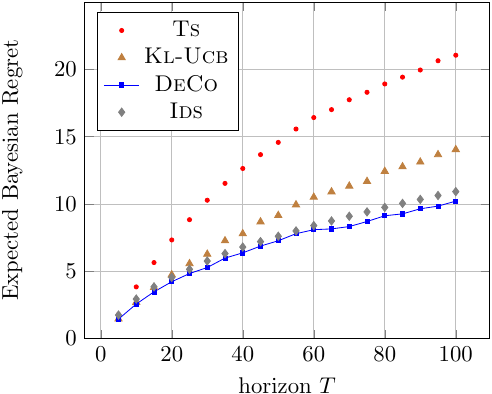}
    \end{subfigure}%
  \end{center}
    \caption{Expected Bayesian regret~\eqref{eq:bayesianRegret} for
      \DeCo\ and a few benchmark policies (the lower the better) for 2, 5, 15 and 20~arms with
      uniform prior.}
    \label{fig:result-main}
\end{figure}

\begin{figure}[hbtp]
  \begin{center}
    \includegraphics[width=0.40\textwidth]{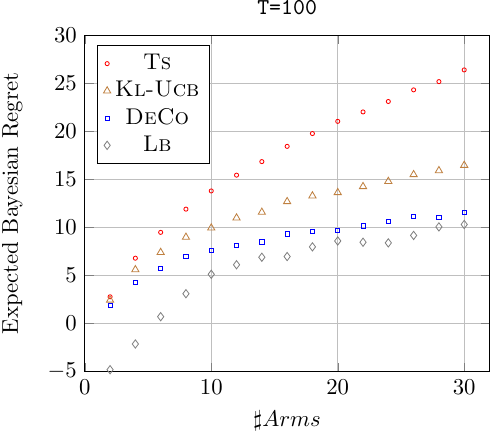}%
    \includegraphics[width=0.40\textwidth]{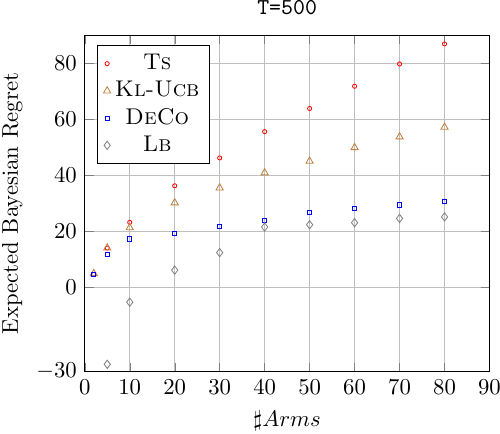}%
  \end{center}
  \caption{Expected Bayesian regret~\eqref{eq:bayesianRegret} for \DeCo, \Ts\
    and \Klucb\, with
    uniform prior, as functions of the number $\sharp {Arms}$ of arms.
    The (\DeCo) lower bound \LowerBound\ in~\eqref{eq:lower_bound_bayesianRegret}
    is also plotted and demonstrates that
    \DeCo\ is close to the optimal solution when the number $\sharp {Arms}$ of
    arms is large enough.
  }
  \label{fig:lowerbound}
\end{figure}

We then tested \DeCo\ 
against Thomson Sampling
(\Ts) \citep{thompson1933likelihood,NIPS2011_e53a0a29},
 Kullback-Leibler upper-confidence bound (\Klucb)
 \citep{garivier2011kl}, Information-Directed Sampling\footnote{For IDS we used the library~\citep{baudry}} 
(\IDS)~\citep{russo2014learning}, Finite Horizon Gittins index~\FH\ \citep{kauffman,nino2011computing,lattimore2016regret}
 and, in the case of two~arms, the exact dynamic programming
 solution\footnote{We did not include Optimistic Gittins
     Index as we did not manage to
     reproduce the results in~\citep{farias2016optimistic}. Also,  it seems that the right hand side  in \cite[example 3.1]{farias2016optimistic} is 
     inexact, as reveals the case $\lambda =1$.}.
Any control sequence
\( \va{\Control}=
\sequence{\va{\Control}_{t}^{\job}}{\job\in\JOB, t\in\ic{0,\horizon-1}} \) is
assessed (and compared to others) using the expected Bayesian regret given by
\begin{align}
  \label{eq:bayesianRegret}
  \REGRET
    \np{ \va{\Control} }=
    \int_{\Delta(\SIMPLEX)^\JOB} \prod_{\job \in \JOB}\prior^{\job}_{0}\np{\dd\proba^{\job}}
    \bigg\{
  \EE_{\sequence{\proba^{\job}}{\job\in\JOB}}
  \bgc{ \sum_{t=0}^{\horizon-1}
    \sum_{\job\in\JOB} \bp{ \va{\Control}_{t}^{\textsc{Ba},\job} - \va{\Control}_{t}^{\job}}
    \va{\Uncertain}_{t+1}^{\job} }\bigg\}
\eqsepv
\end{align}
where we have set the instantaneous costs
$\InstantaneousCost^{\job}_t$ equal to $1$ on \(\good\) and $0$ on \(\bad\) for and where the \textsc{Ba} (best arm) policy is, for all
$\job \in \JOB$, given by
\( \va{\Control}_{t}^{\textsc{Ba},\job}=1 \iff 
\job \in \argmax_{\job'\in \JOB} \proba_{\good}^{\job'} \), 
and where the prior is supposed to be the uniform distribution for all arms.
Numerically, the expected Bayesian regret is obtained by Monte Carlo simulations, where the
expectation with respect to the prior is obtained with a sample of size~$1000$ and expectation with respect to the arms parameters
is obtained with a sample of size~$1000$ in
Figure~\ref{fig:result-main} and of size~$100$ in
Figure~\ref{fig:lowerbound}.
On all cases, \DeCo\ beats both \Ts\ and \Klucb\ with a comfortable
margin, and is comparable to  \IDS. 
For the two arms case in Figure~\ref{fig:result-main}(a), \DeCo\
is very close to the optimal solution, computed by dynamic programming (we used
the Julia BinaryBandit library~\citep{jacko2019finite,jacko2019binarybandit}).

Last, we also computed  the dual bound provided by \DeCo.
Indeed, the upper bound~\eqref{eq:upper_bound_bellman3} yields the inequality
\begin{align}
  \REGRET\np{\va{\Control}}
  \geq
  \REGRET^{\textsc{Lb}}
  =
  \frac{\cardinal{\JOB}}{\cardinal{\JOB}+1}\horizon
  &- \bgp{ \sum_{\job\in\JOB} \VALUE_{0}^{\job}\nc{\multiplier\opt}\np{\prior_{0}^{\job}} 
    +\sum_{t=0}^{\horizon{-}1}\multiplier\opt_t }
  \eqfinp
  \label{eq:lower_bound_bayesianRegret}
\end{align}
Figure~\ref{fig:lowerbound} shows the regret lower
bound, and the \DeCo, \Ts\ and \Klucb\ regrets as a function of the number of
arms for horizons $\horizon=100$
and $\horizon=500$ (beware: in Figure~\ref{fig:lowerbound}, the $x$ axis is the number of arms!).
The lower bound is of no use (lower than 0) for small numbers~2 and 5 of arms.
Nevertheless, when the number of arms increases, the regret of \DeCo\ and the lower bound become quite close,
which indicates that, for those examples, \DeCo\ is close to being optimal.

\begin{table}[hbtp]
  \caption{Comparison of \DeCo, \Bf\ and \FH\ in term of
    estimated total expected reward (higher is better).}
  \label{BFvsDeco}
   \begin{center}
  \begin{tabular}{r|r|r|r}
    \toprule
    $\cardinal{\JOB}$ & $ T$ &     \Bf\ &   \DeCo\ \\
    \midrule
    3 &      10 &   6.409     &     6.411 \\
    3 &      20 &  13.465     &    13.458 \\
    \hline
    5 &      10 &   6.659     &     6.645 \\
    \bottomrule
  \end{tabular}

    \begin{tabular}{r|r|r|r}
    \toprule
    $\cardinal{\JOB}$ & $ T$ &     \FH\ &   \DeCo\  \\
    \midrule
    5 &      20 & 14.28     & 14.21   \\
    5 &     40 &  30.06    &29.85    \\
    \hline
    15 &      20 & 14.67     & 14.59   \\
    15 &     40 &   31.63   & 31.54    \\
    \bottomrule
  \end{tabular}
  \end{center}
\end{table}

\section{Conclusion}

It is notable  that decomposition methods perform so well out of the box,
even when compared to the best known approaches to such a  very scrutinized problem as is the Bayesian bandit. 
The numerical results demonstrate the value of the 
decomposition-coordination approach and shall serve as a proof of
concept: \DeCo\  is a simple algorithm and its performances are close to the optimal
Bayesian solution for several configurations of arms and horizons, while keeping the
computing time reasonable.
Empirically, \DeCo\ offer performances comparable to \FH\,
but with a much smaller computation burden.
For practical applications, we believe \DeCo\ could be useful for decision setting where
the horizon is short and the stack is high.

It should be noted that, as of now, the approach main limitations is that
the horizon~$\horizon$ is supposed to be known in advance and to be  reasonably small (in the experiments, $T\leq 500$), whereas many MABs
algorithms do not require~$\horizon$ as an input.
In addition, the usage of dynamic programming might make \DeCo\ too burdensome
for some applications. 
Also, since  the \DeCo\ algorithm
requires a Bayesian prior, the question of the impact of a wrong prior
on the performance  is left open. 
On the other hand, \DeCo\ can deal with time varying reward functions
and can 
even include a final reward.
In particular, \DeCo\  can be  applied to nonstationary settings
where \FH\ cannot.

On the one hand, this work demonstrates that decomposition methods could be
envisioned for addressing the exploration-exploitation tradeoff. On the other
hand, bandit problems provide an interesting playground to extend the
understanding of those methods, as demonstrated by Sect.~\ref{sec:geometric-representation}.
The characterization of the optimality gap presented in this paper, while focused on a specific problem, is a step toward better understanding  decomposition methods. While it does not fully explain why \DeCo\ works so well, it is a step into that direction as it provides a visual interpretation of the gap between the upper bound provided by the relaxed solution and the performance of \DeCo. 

\newcommand{\noopsort}[1]{} \ifx\undefined\allcaps\def\allcaps#1{#1}\fi

\end{document}